\documentclass[a4paper, 11pt]{article}
\usepackage[T2A,T1]{fontenc}
\usepackage[utf8]{inputenc}
\usepackage[german, russian, english]{babel}
\usepackage{amsmath, amsthm, amssymb,amsfonts}
\usepackage{mathtools}
\usepackage{mathrsfs}
\usepackage{enumerate}
\usepackage{graphicx}
\usepackage{color}
\usepackage{hyperref}
\usepackage[toc, page]{appendix}

\numberwithin{equation}{section}

\newtheorem{theorem}{Theorem}
\newtheorem{lemma}[theorem]{Lemma}
\newtheorem{corollary}[theorem]{Corollary}
\newtheorem{remark}[theorem]{Remark}

\hypersetup{
    pdftitle = {On Stability of Hyperbolic Thermoelastic Reissner-Mindlin-Timoshenko Plates},
    pdfauthor = {Michael Pokojovy},
    pdfsubject = {Thermoelastic plates with hyperbolic heat conduction},
    pdfkeywords = {Reissner-Mindlin-Timoshenko plate; hyperbolic thermoelasticity; second sound; exponential stability; rotational symmetry}
    }

\date{December 4, 2013}

\begin{document}

\title{On Stability of Hyperbolic Thermoelastic Reissner-Mindlin-Timoshenko Plates}
\author{Michael Pokojovy\thanks{E-mail: \texttt{michael.pokojovy@uni-konstanz.de},
Department of Mathematics and Statistics, University of Konstanz, D-78467 Konstanz, Germany}}

\maketitle

\begin{abstract}
	In the present article, we consider a thermoelastic plate of Reissner-Mindlin-Timoshenko type with the hyperbolic heat conduction arising from Cattaneo's law.
	In the absense of any additional mechanical dissipations, the system is often not even strongly stable
	unless restricted to the rotationally symmetric case, etc.
	We present a well-posedness result for the linear problem under general mixed boundary conditions for the elastic and thermal parts.
	For the case of a clamped, thermally isolated plate, we show an exponential energy decay rate under a full damping for all elastic variables.
	Restricting the problem to the rotationally symmetric case, we further prove
	that a single frictional damping merely for the bending compoment is sufficient for exponential stability.
	To this end, we construct a Lyapunov functional incorporating the Bogovski\u i operator for irrotational vector fields
	which we discuss in the appendix.
\end{abstract}

\noindent \textbf{MOS subject classification:} 35L55; 35Q74; 74D05; 93D15; 93D20



\vspace{-6pt}

\section{Introduction}
	Let $\Omega \subset \mathbb{R}^{2}$ be a bounded domain with a boundary $\Gamma := \partial \Omega$.
	We consider a thermoelastic Reissner-Mindlin-Timoshenko plate of a uniform thickness $h > 0$
	such that its midplane occupies the domain $\Omega$ when being in a reference state free of any elastic or thermal stresses.
	The heat propagation is modeled by means of the Cattaneo's law (viz. \cite{Ca1958}).
	With $w$ denoting the vertical displacement of the midplane,
	$\psi$, $\varphi$ the after-bending-angles of vertical filaments being perpendicular to the midplane in the reference state,
	$\theta$ the thermal moment and $q$ the moment of the heat flux, respectively,
	the symmetrized form of Reissner-Mindlin-Timoshenko equations reads as
	\begin{align}
		\rho_{1} w_{tt} - K (w_{x_{1}} + \psi)_{x_{1}} - K (w_{x_{2}} + \varphi)_{x_{2}} &= 0 \text{ in } (0, \infty) \times \Omega,
		\label{GLEICHUNG_REISSNER_MINDLIN_ORIGINAL_LINEAR_1} \\
		\rho_{2} \psi_{tt} - D(\psi_{x_{1} x_{1}} + \tfrac{1 - \mu}{2} \psi_{x_{2} x_{2}} + \tfrac{1 + \mu}{2} \varphi_{x_{1} x_{2}})
		+ K (\psi + w_{x_{1}}) + \gamma \theta_{x_{1}} &= 0 \text{ in } (0, \infty) \times \Omega,
		\label{GLEICHUNG_REISSNER_MINDLIN_ORIGINAL_LINEAR_2} \\
		\rho_{2} \varphi_{tt} - D(\varphi_{x_{2} x_{2}} + \tfrac{1 - \mu}{2} \varphi_{x_{1} x_{1}} + \tfrac{1 + \mu}{2} \psi_{x_{1} x_{2}})
		+ K (\varphi + w_{x_{2}}) + \gamma \theta_{x_{2}} &= 0 \text{ in } (0, \infty) \times \Omega,
		\label{GLEICHUNG_REISSNER_MINDLIN_ORIGINAL_LINEAR_3} \\
		\rho_{3} \theta_{t} + \kappa \mathrm{div}\, q + \beta \theta + \gamma (\psi_{tx_{1}} + \varphi_{tx_{2}}) &= 0 \text{ in } (0, \infty) \times \Omega,
		\label{GLEICHUNG_REISSNER_MINDLIN_ORIGINAL_LINEAR_4} \\
		\tau_{0} q_{t} + \delta q + \kappa \nabla \theta &= 0 \text{ in } (0, \infty) \times \Omega.
		\label{GLEICHUNG_REISSNER_MINDLIN_ORIGINAL_LINEAR_5}
	\end{align}
	A physical deduction of the model can be found in \cite[Kapitel 1]{Po2011}.
	See also \cite[Chapter 1]{LaLi1988} for the case of purely elastic Reissner-Mindlin-Timoshenko plates
	or thermoelastic Kirchhoff-Love plates with parabolic heat conduction.
	Note that, in contrast to the heat equation in uniformly thick bodies,
	$\beta \theta$-term naturally arises in the model
	since $\theta$ is the thermal moment and not the temperature.
	
	Structurally viewed, Reissner-Mindlin-Timoshenko Equations (\ref{GLEICHUNG_REISSNER_MINDLIN_ORIGINAL_LINEAR_1})--(\ref{GLEICHUNG_REISSNER_MINDLIN_ORIGINAL_LINEAR_5})
	can be interpreted as a 2D Lam\'{e} system (\ref{GLEICHUNG_REISSNER_MINDLIN_ORIGINAL_LINEAR_2})--(\ref{GLEICHUNG_REISSNER_MINDLIN_ORIGINAL_LINEAR_3})
	for the filament angles $(\psi, \varphi)'$ coupled to the wave equation (\ref{GLEICHUNG_REISSNER_MINDLIN_ORIGINAL_LINEAR_1}) for the bending component $w$
	and the Cattaneo system (\ref{GLEICHUNG_REISSNER_MINDLIN_ORIGINAL_LINEAR_4})--(\ref{GLEICHUNG_REISSNER_MINDLIN_ORIGINAL_LINEAR_5}) for the thermal moment $\theta$ and the moment of the heat flux $q$.
	Since neither mechanical, no thermal dissipation due to the lack of a direct coupling to the Cattaneo system is present in Equation (\ref{GLEICHUNG_REISSNER_MINDLIN_ORIGINAL_LINEAR_1}),
	one expects the decay properties of (\ref{GLEICHUNG_REISSNER_MINDLIN_ORIGINAL_LINEAR_1})--(\ref{GLEICHUNG_REISSNER_MINDLIN_ORIGINAL_LINEAR_5})
	to be not better than those of classical or hyperbolic 2D thermoelasticity.
	The latter have been investigated by numerous authors.
	Whereas the thermal dissipation arising from the parabolic heat equation leads (with ``few'' exceptions) to the strong stability when coupled with a Lam\'{e} system in a bounded domain of $\mathbb{R}^{n}$ --
	as shown by Dafermos in \cite{Da1968}, no uniform decay can usually be expected (cp. \cite{LeZu1999}).
	Reducing the problem to the case of rotationally symmetric solutions,
	Jiang and Racke \cite[Theorem 4.2]{JiaRa2000} showed an exponential decay of the second-order energy, also in the nonlinear situation (cf. \cite[Theorem 7.3]{JiaRa2000}).
	A similar result was latter obtained by Racke in \cite{Ra2003} for the linear 2D and 3D hyperbolic thermoelasticity.

	As a matter of fact, Reissner-Mindlin-Timoshenko plates and Timoshenko beams have a certain degree of similarity with Kirchhoff-Love plates and Euler-Bernoulli beams.
	The latter also describe the bending of an elastic plate or a beam
	under the assumption that the linear filaments remain perpendicular to mid-plane even after the plate's deformation.
	This model can be shown to be a limit (in a certain sense) of the Timoshenko model as the shear correction factor $K \to \infty$ (cf. \cite{La1989}).
	Numerous mathematical results on Kirchhoff plates are known in the literature.
	In his monograph \cite{La1989}, Lagnese studied various boundary feedback stabilizers furnishing uniform or strong stability for the Kirchhoff-Love plate 
	coupled with a parabolic heat equation in a bounded domain with or without assumptions on the geometry.
	Avalos and Lasiecka exploited further in \cite{AvLa1997} a multiplier technique, interpolation tools and regularity results
	to obtain exponential stability of a thermoelastic Kirchhoff-Love plate without any additional boundary dissipation in the presence or absense of rotational inertia.
	Another important development in this field was made by Lasiecka and Triggiani (see, e.g., \cite{LaTri1989})
	who showed the analyticity of underlying semigroup for all combinations of natural boundary conditions.
	Implying the maximal $L^{2}$-regularity property, this became an important tool for studying nonlinear plates, e.g., the von K\'{a}rm\'{a}n model,
	which was done by Avalos et al. in \cite{AvLaTr1999}.
	It should though be pointed out that this approach is not directly applicable to the case of coupling with the hyperbolic Cattaneo's heat conduction system
	which destroys the analyticity of the semigroup.
	Nonetheless, in an analogous situation of a partly hyperbolic systems such as the full von K\'{a}rm\'{a}n one,
	Lasiecka \cite{La1999} obtained the existence of weak and regular solutions
	and showed their uniform stability in the presense of a mechanical damping only for the solenoidal part for the in-plane displacements.
	A similar study has then later been carried out in \cite{BeLa2000} by Benabdallah and Lasiecka for the full von K\'{a}rm\'{a}n model incorporating rotational inertia.

	Turning back to Reissner-Mindlin-Timoshenko plates, we once again refer to the monograph \cite{La1989} of Lagnese 
	in which he addressed the question of uniform (in particular, exponential) and strong stabilization of purely elastic plates by the means of boundary feedbacks.
	For the following choice of stabilizing feedbacks on a portion $\Gamma_{1} \neq \emptyset$ of the boundary
	\begin{align}
		\begin{split}
			w = \psi = \varphi &= 0 \;\;\, \text{ in } (0, \infty) \times \Gamma_{0}, \\
			K(\tfrac{\partial w}{\partial \nu} + \nu_{1} \psi + \nu_{2} \varphi) &= m_{1} \text{ in } (0, \infty) \times \Gamma_{1}, \\
			D(\nu_{1} \psi_{x_{1}} + \mu \nu_{1} \varphi_{x_{2}} + \tfrac{1 - \mu}{2} (\psi_{x_{2}} + \varphi_{x_{1}}) \nu_{2}) &= m_{2} \text{ in } (0, \infty) \times \Gamma_{1}, \\
			D(\nu_{2} \varphi_{x_{2}} + \mu \nu_{2} \psi_{x_{1}} + \tfrac{1 - \mu}{2} (\psi_{x_{2}} + \varphi_{x_{1}}) \nu_{1}) &= m_{3} \text{ in } (0, \infty) \times \Gamma_{1},
		\end{split} \notag
	\end{align}
	the purely elastic Reissner-Mindlin-Timoshenko plate
	\begin{align}
		\begin{split}
			\rho_{1} w_{tt} - K (w_{x_{1}} + \psi)_{x_{1}} - K (w_{x_{2}} + \varphi)_{x_{2}} &= 0 \text{ in } (0, \infty) \times \Omega, \\
			\rho_{2} \psi_{tt} - D(\psi_{x_{1} x_{1}} + \tfrac{1 - \mu}{2} \psi_{x_{2} x_{2}} + \tfrac{1 + \mu}{2} \varphi_{x_{1} x_{2}})
			+ K (\psi + w_{x_{1}}) + \gamma \theta_{x_{1}} + d_{1} \psi_{t} &= 0 \text{ in } (0, \infty) \times \Omega, \\
			\rho_{2} \varphi_{tt} - D(\varphi_{x_{2} x_{2}} + \tfrac{1 - \mu}{2} \varphi_{x_{1} x_{1}} + \tfrac{1 + \mu}{2} \psi_{x_{1} x_{2}})
			+ K (\varphi + w_{x_{2}}) + \gamma \theta_{x_{2}} + d_{2} \varphi_{t} &= 0 \text{ in } (0, \infty) \times \Omega
		\end{split}
		\notag
	\end{align}
	was proved to be strongly stable (i.e., the energy was shown to vanish as $t \to \infty$)
	if $\Gamma_{0} \neq \varnothing$ und $(m_{1}, m_{2}, m_{3})' = -F (w_{t}, \psi_{t}, \varphi_{t})'$
	where $F \in L^{\infty}(\Gamma_{1}, \mathbb{R}^{3 \times 3})$ is a symmetric positive semidefinite matrix function
	which is additionally positive definite on a connected nontrivial portion of $\Gamma_{1}$, etc.
	Under the geometric condition stating that $(\Omega, \Gamma_{0}, \Gamma_{1})$ is ``star complemented---star shaped''
	and some additional assumptions on $F$,
	even uniform stability has been shown.

	Similar results were also obtained by Mu\~{n}oz Rivera and Portillo Oquendo in \cite{MuRiPoOq2003} auch for the boundary conditions of memory-type
	\begin{align}
		\begin{split}
			w = \psi = \varphi &= 0 \text{ in } (0, \infty) \times \Gamma_{0}, \\
			w + \int_{0}^{t} g_{1}(t - s) K(\tfrac{\partial w}{\partial \nu} + \nu_{1} \psi + \nu_{2} \varphi)(s) \mathrm{d} s &= 0 \text{ in } (0, \infty) \times  \Gamma_{1}, \\
			\psi + \int_{0}^{t} g_{1}(t - s) D(\nu_{1} \psi_{x_{1}} + \mu \nu_{1} \varphi_{x_{2}} + \tfrac{1 - \mu}{2} (\psi_{x_{2}} + \varphi_{x_{1}})(s) \mathrm{d} s &= 0 \text{ in } (0, \infty) \times \Gamma_{1}, \\
			\varphi + \int_{0}^{t} g_{1}(t - s) D(\nu_{2} \varphi_{x_{2}} + \mu \nu_{2} \psi_{x_{1}} + \tfrac{1 - \mu}{2} (\psi_{x_{2}} + \varphi_{x_{1}}) \nu_{1})(s) \mathrm{d} s &= 0 \text{ in } (0, \infty) \times \Gamma_{1}
		\end{split} \notag
	\end{align}
	with exponential kernels $g_{1}$, $g_{2}$, $g_{3}$.

	In \cite{FeSa2009}, Fern\'{a}ndez Sare studied a linear Reissner-Mindlin-Timoshenko plate with a damping for both angle components
	\begin{align}
		\rho_{1} w_{tt} - K (w_{x_{1}} + \psi)_{x_{1}} - K (w_{x_{2}} + \varphi)_{x_{2}} &= 0 \text{ in } (0, \infty) \times \Omega,
		\label{GLEICHUNG_REISNNER_MINDLIN_ELASTISCH_1} \\
		\rho_{2} \psi_{tt} - D(\psi_{x_{1} x_{1}} + \tfrac{1 - \mu}{2} \psi_{x_{2} x_{2}} + \tfrac{1 + \mu}{2} \varphi_{x_{1} x_{2}})
		+ K (\psi + w_{x_{1}}) + \gamma \theta_{x_{1}} + d_{1} \psi_{t} &= 0 \text{ in } (0, \infty) \times \Omega,
		\label{GLEICHUNG_REISNNER_MINDLIN_ELASTISCH_2} \\
		\rho_{2} \varphi_{tt} - D(\varphi_{x_{2} x_{2}} + \tfrac{1 - \mu}{2} \varphi_{x_{1} x_{1}} + \tfrac{1 + \mu}{2} \psi_{x_{1} x_{2}})
		+ K (\varphi + w_{x_{2}}) + \gamma \theta_{x_{2}} + d_{2} \varphi_{t} &= 0 \text{ in } (0, \infty) \times \Omega.
		\label{GLEICHUNG_REISNNER_MINDLIN_ELASTISCH_3}
	\end{align}
	He proved that the system is polynomially stable under Dirichlet boundary conditions on all three variables.
	For a particular choice of boundary conditions in a rectangular configuration $\Omega = (0, L_{1}) \times (0, L_{2})$,
	a resolvent criterion was exploited to show that the system is not exponentially stable.
	
	Mu\~{n}oz Rivera und Racke considered in \cite{MuRiRa2002} an nonlinear Timoshenko-beam coupled to a parabolic heat equation
	\begin{equation}
		\begin{split}
			\rho_{1} \varphi_{tt} - \sigma(\varphi_{x}, \psi)_{x} &= 0 \text{ in } (0, \infty) \times (0, L), \\
			\rho_{2} \psi_{tt} - b \psi_{xx} + k(\varphi_{x} + \psi) + \gamma \theta_{x} &= 0 \text{ in } (0, \infty) \times (0, L), \\
			\rho_{3} \theta_{t} - \kappa \theta_{xx} + \gamma \psi_{tx} &= 0 \text{ in } (0, \infty) \times (0, L) \\
		\end{split} \notag
	\end{equation}
	subject to mixed boundary conditions
	$\varphi = \psi = \theta_{x} = 0$ or $\varphi = \psi_{x} = \theta = 0$.
	Both in the linear case, i.e., $\sigma(r, s) = kr + s$,
	and the nonlinear case, i.e., for a smooth stress function $\sigma$ satisfying $\nabla \sigma = (k, k)'$, $\nabla^{2} \sigma = 0$, but in the latter case only for sufficiently small initial data,
	the energy was shown to decay exponentially if the condition $\tfrac{\rho_{1}}{k} = \tfrac{\rho_{2}}{b}$ holds true.
	For the linear situation, this condition was even shown to be necessary for the exponential stability.
	It should though be pointed out that the latter proportionality condition, being mathematically fully sound, is physically not possible.
	
	Surprisingly, this result could not be carried over to the case of Cattaneo heat conduction.
	Namely, Fern\'{a}ndez Sare and Racke showed in \cite{FeSaRa2009} that the purely hyperbolic system
	\begin{align}
		\rho_{1} \varphi_{tt} - k(\varphi_{x} + \psi)_{x} &= 0 \text{ in } (0, \infty) \times (0, L),
		\label{GLEICHUNG_TIMOSHENKO_1} \\
		\rho_{2} \psi_{tt} - b \psi_{xx} + k(\varphi_{x} + \psi) + \gamma \theta_{x} &= 0 \text{ in } (0, \infty) \times (0, L),
		\label{GLEICHUNG_TIMOSHENKO_2} \\
		\rho_{3} \theta_{t} + \kappa q_{x} + \gamma \psi_{tx} &= 0 \text{ in } (0, \infty) \times (0, L),
		\label{GLEICHUNG_TIMOSHENKO_3} \\
		\tau_{0} q_{t} + \delta q + \kappa \theta_{x} &= 0  \text{ in } (0, \infty) \times (0, L)
		\label{GLEICHUNG_TIMOSHENKO_4}
	\end{align}
	is not exponentially stable even under the assumption $\tfrac{\rho_{1}}{k} = \tfrac{\rho_{2}}{b}$.
	This motivated Messaoudi et al. to introduce a frictional damping for the bending component.
	In \cite{MePoSa2009}, they replaced Equation (\ref{GLEICHUNG_TIMOSHENKO_1}) with the damped equation
	\begin{equation}
		\rho_{1} \varphi_{tt} - \sigma(\varphi_{x}, \psi)_{x} + \mu \varphi_{t} = 0 \text{ in } (0, \infty) \times (0, L) \notag
	\end{equation}
	for some $\mu > 0$.
	Under this additional mechanical dissipation, they proved that both linear and nonlinear systems are stable
	under the boundary conditions $\varphi = \psi = q = 0$ und $\varphi_{x} = \psi = q = 0$
	independent of whether the relation $\tfrac{\rho_{1}}{k} = \tfrac{\rho_{2}}{b}$ holds or not.

	The impact of thermal coupling on the strong stability of a Reissner-Mindlin-Timoshenko plate has also been studied by Grobbelaar in her papers \cite{Gro2011}, \cite{Gro2012} and \cite{Gro2013}.
	In \cite{Gro2011}, the author considered a stuctural 3D acoustic model with a 2D plate interface
	and proved a strong asymptotic stability for the radially symmetric case.
	A similar result was later obtained in \cite{Gro2012} for a rotationally symmetric Reissner-Mindlin-Timoshenko plate with hyperbolic heat conduction due to Cattaneo.
	To this end, both articles employed Benchimol's spectral criterion. The arguments can be directly carried over to the case of classical Fourier heat conduction
	being a formal limit Cattaneo's system as the relaxation parameter $\tau \to 0$.
	In her recent article \cite{Gro2013}, Grobbelaar proved a polynimal decay rate of $t^{-1/4}$
	in the rotationally symmetric case for the Reissner-Mindlin-Timoshenko system coupled to the classical Fourier heat conduction
	under Dirichlet boundary conditions on $w$ and $\theta$ as well as free boundary conditions on $\psi$ and $\varphi$.

	In the present article, we consider the linear Reissner-Mindlin-Timosheko plate equations (\ref{GLEICHUNG_REISSNER_MINDLIN_ORIGINAL_LINEAR_1})--(\ref{GLEICHUNG_REISSNER_MINDLIN_ORIGINAL_LINEAR_5}) in a bounded domain.
	The paper is structured as follows.
	In the first section, we exploit the semigroup theory to show that the initial-boundary value problem (\ref{GLEICHUNG_REISSNER_MINDLIN_ORIGINAL_LINEAR_1})--(\ref{GLEICHUNG_REISSNER_MINDLIN_ORIGINAL_LINEAR_5})
	subject to corresponding initial conditions as well as homogeneous Dirichlet and Neumann boundary conditions on both elastic and thermal variables on different portions of the boundary is well-posed.
	In the second section, we prove the lack of strong stability for this problem provided $\Gamma$ is smooth for a particular set of boundary conditions.
	We further show that a mechanical damping for all three variables $w$, $\varphi$ and $\psi$ leads to an exponential decay rate
	under Dirichlet boundary conditions for the elastic and Neumann boundary conditions for the thermal part of the system.
	Restricting the domain $\Omega$ and the data to the rotationally symmetric case,
	we prove that a single mechanical damping on $w$ is enough to exponentially stabilize the system.
	This is a generalization of Messaoudi's et al. stability results from \cite{MePoSa2009} to a multi-dimensional situation.
	In the appendix, we finally present a brief discussion on Bogovski\u i operator for irrotational vector fields and show its continuity.

\section{Existence and uniqueness of classical solutions}
	In the following, unless specified otherwise, we assume the boundary $\Gamma$ to be Lipschitzian and satisfy
	$\Gamma = \bar{\Gamma}_{1} \cup \bar{\Gamma}_{2} = \bar{\Gamma}_{3} \cup \bar{\Gamma}_{4}$
	with $\Gamma_{1} \neq \varnothing$, $\Gamma_{1} \cap \Gamma_{2} = \emptyset$, $\Gamma_{3} \cap \Gamma_{4} = \emptyset$ and $\Gamma_{k}$, $k = 1, \dots, 4$, being relatively open.
	Let the plate be clamped at $\Gamma_{1}$ and hinged at $\Gamma_{2}$.
	Further, let it be held at the reference temperature on $\Gamma_{3}$
	and be thermally insulated on $\Gamma_{4}$.
	Then, the boundary conditions read as
	\begin{align}
		w = \psi = \varphi &= 0 \text{ on } (0, \infty) \times \Gamma_{1},
		\label{GLEICHUNG_REISSNER_MINDLIN_LINEAR_RB_1} \\
		K(\tfrac{\partial w}{\partial \nu} + \nu_{1} \psi + \nu_{2} \varphi) &= 0 \text{ on } (0, \infty) \times \Gamma_{2},
		\label{GLEICHUNG_REISSNER_MINDLIN_LINEAR_RB_2} \\
		D(\nu_{1} \psi_{x_{1}} + \mu \nu_{1} \varphi_{x_{2}} + \tfrac{1 - \mu}{2} (\psi_{x_{2}} + \varphi_{x_{1}}) \nu_{2}) - \gamma \theta \nu_{1} &= 0 \text{ on } (0, \infty) \times \Gamma_{2},
		\label{GLEICHUNG_REISSNER_MINDLIN_LINEAR_RB_3} \\
		D(\nu_{2} \varphi_{x_{2}} + \mu \nu_{2} \psi_{x_{1}} + \tfrac{1 - \mu}{2} (\psi_{x_{2}} + \varphi_{x_{1}}) \nu_{1}) - \gamma \theta \nu_{2} &= 0 \text{ on } (0, \infty) \times \Gamma_{2},
		\label{GLEICHUNG_REISSNER_MINDLIN_LINEAR_RB_4} \\
		\theta &= 0 \text{ on } (0, \infty) \times \Gamma_{3},
		\label{GLEICHUNG_REISSNER_MINDLIN_LINEAR_RB_5} \\
		q \cdot \nu &= 0 \text{ on } (0, \infty) \times \Gamma_{4},
		\label{GLEICHUNG_REISSNER_MINDLIN_LINEAR_RB_6}
	\end{align}
	where $\nu = (\nu_{1}, \nu_{2})'$ denotes the outer unit normal vector to $\Gamma$
	and $(\cdot)'$ stands for the usial matrix transposition.
	
	Using the standard notation from the Theory of elasticity (cf. \cite[p. 8]{JiaRa2000}),
	we introduce the generalized gradient and the corresponding boundary symbol
	\begin{align}
		\mathcal{D} := \left(\begin{array}{cc}
			\partial_{1} & 0 \\
			0 & \partial_{2} \\
			\partial_{2} & \partial_{1}
		\end{array}\right), \quad
		\mathcal{N} := \left(\begin{array}{cc}
			\nu_{1} & 0 \\
			0 & \nu_{2} \\
			\nu_{2} & \nu_{1}
		\end{array}\right), \notag
	\end{align}
	respectively.
	With this notation, we can easily conclude
	\begin{equation}
		\begin{split}
			D \left(\begin{array}{c}
				\psi_{x_{1} x_{1}} + \tfrac{1 - \mu}{2} \psi_{x_{2} x_{2}} + \tfrac{1 + \mu}{2} \varphi_{x_{1} x_{2}} \\
				\varphi_{x_{2} x_{2}} + \tfrac{1 - \mu}{2} \varphi_{x_{1} x_{2}} + \tfrac{1 + \mu}{2} \psi_{x_{1} x_{2}}
			\end{array}\right)
			&= \mathcal{D}' S \mathcal{D} v, \\
			D \left(\begin{array}{c}
				\nu_{1} \psi_{x_{1}} + \mu \nu_{1} \varphi_{x_{2}} + \tfrac{1 - \mu}{2} (\psi_{x_{2}} + \varphi_{x_{1}}) \nu_{2} \\
				\nu_{2} \varphi_{x_{2}} + \mu \nu_{2} \psi_{x_{1}} + \tfrac{1 + \mu}{2} (\psi_{x_{2}} + \varphi_{x_{1}}) \nu_{1}
			\end{array}\right)
			&= \mathcal{N}' S \mathcal{D} v,
		\end{split}
		\notag
	\end{equation}
	where $v := (\psi, \varphi)'$ and
	\begin{equation}
		S := D \left(\begin{array}{ccc}
			1 & \mu & 0 \\
			\mu & 1 & 0 \\
			0 & 0 & \tfrac{1 - \mu}{2}
		\end{array}\right).
		\label{GLEICHUNG_STEIFIGKEITSMATRIX}
	\end{equation}
	With $\mu$ satisfying $\mu \in (-1, 1)$,
	the symmetric matrix $S$ is positive definite since
	$\sigma(S) = \{D \tfrac{1 - \mu}{2}, D (1 - \mu), D (1 + \mu)\}$ due to the fact
	\begin{equation}
		\det(S - \lambda I) = (D \, \tfrac{1 - \mu}{2} - \lambda) ((D - \lambda)^{2} - \mu^{2} D^{2})
		= (D \, \tfrac{1 - \mu}{2} - \lambda) (D - \lambda - \mu D) (D - \lambda + \mu D). \notag
	\end{equation}
	By the virtue of physical condition $\mu \in (0, \tfrac{1}{2})$,
	the latter is not an actual restriction.
	Throughout this section, we assume $S$ to be an arbitrary symmetric, positive definite matrix, i.e., $S \in \mathrm{SPD}(\mathbb{R}^{3})$.
	
	With the notations above,
	Equations (\ref{GLEICHUNG_REISSNER_MINDLIN_ORIGINAL_LINEAR_1})--(\ref{GLEICHUNG_REISSNER_MINDLIN_ORIGINAL_LINEAR_5})
	can be equivalently written as
	\begin{align}
		\rho_{1} w_{tt} - K \mathrm{div}\, (\nabla w + v) &= 0 \text{ in } (0, \infty) \times \Omega, 
		\label{GLEICHUNG_REISNNER_MINDLIN_LINEAR_KOMPAKTE_FORM_1} \\
		\rho_{2} v_{tt} - \mathcal{D}' S\mathcal{D} v + K(v + \nabla w) + \gamma \nabla \theta &= 0 \text{ in } (0, \infty) \times \Omega,
		\label{GLEICHUNG_REISNNER_MINDLIN_LINEAR_KOMPAKTE_FORM_2} \\
		\rho_{3} \theta_{t} + \kappa \mathrm{div}\, q + \beta \theta + \gamma \mathrm{div}\, v_{t} &= 0 \text{ in } (0, \infty) \times \Omega,
		\label{GLEICHUNG_REISNNER_MINDLIN_LINEAR_KOMPAKTE_FORM_3} \\
		\tau_{0} q_{t} + \delta q + \kappa \nabla \theta  &= 0 \text{ in } (0, \infty) \times \Omega
		\label{GLEICHUNG_REISNNER_MINDLIN_LINEAR_KOMPAKTE_FORM_4}
	\end{align}
	with the boundary conditions (\ref{GLEICHUNG_REISSNER_MINDLIN_LINEAR_RB_1})--(\ref{GLEICHUNG_REISSNER_MINDLIN_LINEAR_RB_6}) transformed to
	\begin{align}
		w = |v| &= 0 \text{ on } (0, \infty) \times \Gamma_{1}, 
		\label{GLEICHUNG_REISNNER_MINDLIN_LINEAR_KOMPAKTE_FORM_RB_1} \\
		(\nabla w + v) \cdot \nu &= 0 \text{ on } (0, \infty) \times \Gamma_{2},
		\label{GLEICHUNG_REISNNER_MINDLIN_LINEAR_KOMPAKTE_FORM_RB_2} \\
		\mathcal{N}' S \mathcal{D} v - \gamma \theta \nu &= 0 \text{ on } (0, \infty) \times \Gamma_{2},
		\label{GLEICHUNG_REISNNER_MINDLIN_LINEAR_KOMPAKTE_FORM_RB_3} \\
		\theta &= 0 \text{ on } (0, \infty) \times \Gamma_{3},
		\label{GLEICHUNG_REISNNER_MINDLIN_LINEAR_KOMPAKTE_FORM_RB_4} \\
		q \cdot \nu &= 0 \text{ on } (0, \infty) \times \Gamma_{4}
		\label{GLEICHUNG_REISNNER_MINDLIN_LINEAR_KOMPAKTE_FORM_RB_5}
	\end{align}
	and initial conditions
	\begin{equation}
		w(0, \cdot) = w^{0}, \; w_{t}(0, \cdot) = w^{1}, \;
		v(0, \cdot) = v^{0}, \; v_{t}(0, \cdot) = v^{1}, \;
		\theta(0, \cdot) = \theta^{0}, \; q(0, \cdot) = q^{0},
		\label{GLEICHUNG_REISNNER_MINDLIN_LINEAR_KOMPAKTE_FORM_AB}
	\end{equation}
	where $v^{0} = (\psi^{0}, \varphi^{0})'$, $v^{1} = (\psi^{1}, \varphi^{1})'$.

\subsection{Well-Posedness}
	\label{ABSCHNITT_REISSNER_MINDLIN_WOHLGESTELLTHEIT}
	We further exploit the semigroup theory to obtain the classical well-posedness of Reissner-Mindlin-Timoshenko equations.
	To this end, we transform Equations (\ref{GLEICHUNG_REISNNER_MINDLIN_LINEAR_KOMPAKTE_FORM_1})--(\ref{GLEICHUNG_REISNNER_MINDLIN_LINEAR_KOMPAKTE_FORM_AB})
	into the Cauchy problem
	\begin{equation}
		\begin{split}
			\frac{\mathrm{d}}{\mathrm{d}t} V(t) &= \mathcal{A} V(t) \text{ f\"ur } t \in (0, \infty), \\
			V(0) &= V^{0}
		\end{split} \notag
	\end{equation}
	on a Hilbert space $\mathcal{H}$.
	According to \cite[Theorem 1.3]{Pa1992},
	the latter is well-posed if and only if
	$\mathcal{A}$ is an infinitesimal generator of a strongly continuous semigroup on $\mathcal{H}$.

	We set $V := (w, v, w_{t}, v_{t}, \theta, q)'$ and formally define the differential operator
	\begin{equation}
		A := \rho^{-1}
		\left(\begin{array}{cccccccc}
			0 & 0 & 1 & 0 & 0 & 0 \\
			0 & 0 & 0 & 1 & 0 & 0 \\
			K \triangle & K \mathrm{div}\, & 0 & 0 & 0 & 0 \\
			-K \nabla & \mathcal{D}'S\mathcal{D} - K & 0 & 0 & -\gamma \nabla & 0 \\
			0 & 0 & 0 & -\gamma \mathrm{div}\, & -\beta & -\kappa \mathrm{div}\, \\
			0 & 0 & 0 & 0 & -\kappa \nabla & -\delta
		\end{array}\right) \notag
	\end{equation}
	with $\rho := \mathrm{diag}(1, 1, \rho_{1}, \rho_{2}, \rho_{3}, \tau_{0})$.
	To introduce the functional analytic settings, we consider the Hilbert space
	\begin{equation}
		\mathcal{H} := (H^{1}_{\Gamma_{1}}(\Omega))^{3} \times (L^{2}(\Omega))^{3} \times (L^{2}(\Omega))^{3} \notag
	\end{equation}
	equipped with the scalar product
	\begin{equation}
		\begin{split}
			\langle V, W\rangle_{\mathcal{H}} := \;&
			\rho_{1} \langle V^{3}, W^{3}\rangle_{L^{2}(\Omega)} + \rho_{2} \langle V^{4}, W^{4}\rangle_{(L^{2}(\Omega))^{2}} +
			K \langle \nabla V^{1} + V^{2}, \nabla W^{1} + W^{2}\rangle_{(L^{2}(\Omega))^{2}} + \\
			&\langle \mathcal{D} V^{2}, S \mathcal{D} W^{2}\rangle_{(L^{2}(\Omega))^{3}} +
			\rho_{3} \langle V^{5}, W^{5}\rangle_{(L^{2}(\Omega))^{2}} + \tau_{0} \langle V^{6}, W^{6}\rangle_{L^{2}(\Omega)}.
		\end{split} \notag
	\end{equation}
	Here, we define for a relatively open set $\Gamma_{0} \subset \Gamma$
	\begin{equation}
		H^{1}_{\Gamma_{0}}(\Omega)
		= \mathrm{cl}\left(\{u \in \mathcal{C}^{\infty}(\Omega) \,|\, \mathrm{supp}(u) \cap \Gamma_{0} = \emptyset\}, \|\cdot\|_{H^{1}(\Omega)}\right). \notag
	\end{equation}
	Note that due to the Lipschitz continuity of $\Gamma$,
	there exists a linear, continuous operator $T \colon H^{1}(\Omega) \to H^{1/2}(\Gamma)$.
	Thus, the notation $u|_{\Gamma_{0}} = 0$ is also legitimate.

	The following theorem implies that $\langle \cdot, \cdot\rangle_{\mathcal{H}}$ is equivalent with the standard product topology on $\mathcal{H}$, i.e., $\mathcal{H}$ is complete.
	The proof is a direct consequence of an analogous result in \cite{LaLi1988}
	(cf. also \cite{Lei1986} for the case of domains with a strict cone property).
	\begin{lemma}
		\label{SATZ_KORNSCHE_UNGLEICHUNG}
		There exist constants $C_{\mathcal{K}, 1}, C_{\mathcal{K}, 2}, C_{\mathcal{K}} > 0$ such that
		\begin{equation}
			C_{\mathcal{K}, 1} \|v\|_{(H^{1}(\Omega))^{2}} \leq \|\sqrt{S} \mathcal{D} v\|_{(L^{2}(\Omega))^{3}} \leq C_{\mathcal{K}, 2} \|v\|_{(H^{1}(\Omega))^{2}} \notag
		\end{equation}
		and
		\begin{equation}
			\|\sqrt{S} \mathcal{D} v\|^{3}_{(L^{2}(\Omega))^{2}} + K \|\nabla w + v\|^{2}_{(L^{2}(\Omega))^{2}}
			\geq C_{\mathcal{K}} \big(\|v\|^{2}_{(H^{1}(\Omega))^{2}} + \|w\|^{2}_{H^{1}(\Omega)}\big) \notag
		\end{equation}
		holds for any  $(w, v) \in (H^{1}_{\Gamma_{1}}(\Omega))^{3}$. 
	\end{lemma}
	
	We introduce the operator
	\begin{equation}
		\mathcal{A} \colon D(\mathcal{A}) \subset \mathcal{H} \longrightarrow \mathcal{H}, \quad V \longmapsto AV, \notag
	\end{equation}
	where
	\begin{align}
		D(\mathcal{A}) = \{V \in \mathcal{H} \,|\, &A V \in \mathcal{H}, V \text{ satisfies the generalized Neumann boundary conditions (\ref{GLEICHUNG_NEUMANN_MR_LIN_KF_1})--(\ref{GLEICHUNG_NEUMANN_MR_LIN_KF_3})}\} \notag \\
		= \{V \in \mathcal{H} \,|\, &V^{1}, V^{3} \in H^{1}_{\Gamma_{1}(\Omega)}, V^{2}, V^{4} \in (H^{1}_{\Gamma_{1}(\Omega)})^{2},
		\triangle V^{1} \in L^{2}(\Omega),
		\mathcal{D}^{T} S \mathcal{D} V^{2} \in (L^{2}(\Omega))^{2}, \notag \\
		&V^{5} \in H^{1}_{\Gamma_{3}}(\Omega), \mathrm{div}\, V^{6} \in L^{2}(\Omega), \notag \\
		&V \text{ satisfies the generalized Neumann boundary conditions (\ref{GLEICHUNG_NEUMANN_MR_LIN_KF_1})--(\ref{GLEICHUNG_NEUMANN_MR_LIN_KF_3})}\} \notag
	\end{align}
	with the generalized Neumann boundary conditions given by
	\begin{align}
		\langle \triangle V^{1} + \mathrm{div}\, V^{2}, \phi\rangle_{L^{2}(\Omega)} + \langle \nabla V^{1} + V^{2}, \nabla \phi \rangle_{(L^{2}(\Omega))^{2}}
		&= 0 \text{ for all } \phi \in H^{1}_{\Gamma_{1}}(\Omega) 
		\label{GLEICHUNG_NEUMANN_MR_LIN_KF_1} \\
		\langle \mathcal{D}^{T} S \mathcal{D} V^{2} - \gamma \nabla V^{5}, \phi\rangle_{(L^{2}(\Omega))^{2}}
		+ \langle S \mathcal{D} V^{2}, \mathcal{D} \phi \rangle_{(L^{2}(\Omega))^{3}} & \notag \\
		- \gamma \langle V^{5}, \mathrm{div}\, \phi \rangle_{L^{2}(\Omega)}
		&= 0 \text{ for all } \phi \in (H^{1}_{\Gamma_{1}}(\Omega))^{2}
		\label{GLEICHUNG_NEUMANN_MR_LIN_KF_2} \\
		\langle \mathrm{div}\, V^{6}, \phi\rangle_{L^{2}(\Omega)} + \langle V^{6}, \nabla \phi\rangle_{(L^{2}(\Omega))^{2}} &= 0 \text{ for all } \phi \in H^{1}_{\Gamma_{3}}(\Omega).
		\label{GLEICHUNG_NEUMANN_MR_LIN_KF_3}
	\end{align}
	Obviously, $D(\mathcal{A})$ is a linear subspace of $\mathcal{H}$.
	
	Thus,
	\begin{equation}
		V_{t} = \mathcal{A} V, \quad V(0) = V^{0} \label{REISSNER_MINDLIN_EVOLUTIONSPROBLEM}
	\end{equation}
	is a generalization of (\ref{GLEICHUNG_REISNNER_MINDLIN_LINEAR_KOMPAKTE_FORM_1})--(\ref{GLEICHUNG_REISNNER_MINDLIN_LINEAR_KOMPAKTE_FORM_AB})
	since any classically differentiable solution to (\ref{GLEICHUNG_REISNNER_MINDLIN_LINEAR_KOMPAKTE_FORM_1})--(\ref{GLEICHUNG_REISNNER_MINDLIN_LINEAR_KOMPAKTE_FORM_AB})
	solves the abstract Cauchy problem (\ref{REISSNER_MINDLIN_EVOLUTIONSPROBLEM}).
	Here, $V^{0} := (w^{0}, v^{0}, w^{1}, v^{1}, \theta^{0}, q^{0})'$ is assumed to be an element of $D(\mathcal{A})$.

	The following theorem characterizes $\mathcal{A}$ as an infinitesimal generator of a strongly continuous semigroup of bounded linear operators on $\mathcal{H}$.
	\begin{theorem}
		The following statements hold true for $\mathcal{A}$.
		\begin{enumerate}
			\item $D(\mathcal{A})$ is dense in $\mathcal{H}$.
			\item $\mathcal{A}$ is a closed operator.
			\item $\mathrm{im}(\lambda - \mathcal{A}) = \mathcal{H}$ for any $\lambda > 0$.
			\item $\mathcal{A}$ is dissipative.
		\end{enumerate}
	\end{theorem}

	\begin{proof}
		\begin{enumerate}
			\item The fact that $D(\mathcal{A})$ is a dense subspace of $\mathcal{H}$ is a direct consequence of the inclusion
			\begin{equation}
				(\mathcal{C}^{\infty}(\Omega))^{9} \cap \mathcal{H} \subset D(\mathcal{A}). \notag
			\end{equation}
			Note that the generalized Neumann boundary conditions (\ref{GLEICHUNG_NEUMANN_MR_LIN_KF_1})--(\ref{GLEICHUNG_NEUMANN_MR_LIN_KF_3}) are satisfied per defintion.
	
			\item The proof of the closedness of $\mathcal{A}$ is also standard.
			We select an arbitrary sequence
			$(V_{n})_{n \in \mathbb{N}} \subset D(\mathcal{A})$ such that $V_{n} \to V \in \mathcal{H}$ and $\mathcal{A} V_{n} \to F \in \mathcal{H}$ as $n \to \infty$
			and show that $V \in D(\mathcal{A})$ and $\mathcal{A} V = F$
			(cf. \cite{Po2011} for the case $\Gamma_{2} = \Gamma_{3} = \emptyset$).
			
			Taking into account $((L^{2}(\Omega))^{9})' \subset \mathcal{H}'$,
			the strong convergence in $\mathcal{H}$ implies the weak convergence in $(L^{2}(\Omega))^{9}$, i.e.,
			\begin{equation}
				\langle \mathcal{A} V_{n}, \Phi\rangle_{(L^{2}(\Omega))^{9}} \to \langle F, \Phi\rangle_{(L^{2}(\Omega))^{9}} \text{ as } n \to \infty \notag
			\end{equation}
			for any $\Phi \in (L^{2}(\Omega))^{9}$.
			With a proper selection of $\Phi$, the problem can be projected onto a corresponding component.
			The proof will be made by means of a proper selection of $\Phi$.
		
			There generally holds for $V \in D(\mathcal{A})$
			\begin{equation}
				\mathcal{A} V = \rho^{-1} \left(\begin{array}{c}
					V^{3} \\
					V^{4} \\
					K \triangle V^{1} + K \mathrm{div}\, V^{2} \\
					-K \nabla V^{1} + \mathcal{D}' S \mathcal{D} V^{2} - K V^{2} - \gamma \nabla V^{5} \\
					-\gamma \mathrm{div}\, V^{4} - \beta V^{5} - \kappa \mathrm{div}\, V^{6} \\
					-\kappa \nabla V^{5} - \delta V^{6}
				\end{array}\right). \notag
			\end{equation}
		
			We consider the following cases:
			\begin{itemize}
				\newcounter{zaehler}
				\setcounter{zaehler}{1}
				\item[{\it \roman{zaehler})}] \addtocounter{zaehler}{1}
				First, we select $\Phi = (\phi, 0, 0, 0, 0, 0)'$, $\phi \in H^{1}_{\Gamma_{1}}(\Omega)$ to obtain
				\begin{align}
					\langle F^{1}, \phi\rangle_{L^{2}(\Omega)} = \langle F, \Phi\rangle_{(L^{2}(\Omega))^{9}}
					\leftarrow \langle AV_{n}, \Phi\rangle_{(L^{2}(\Omega))^{9}} =
					\tfrac{1}{\rho_{1}} \langle V^{3}_{n}, \phi\rangle_{L^{2}(\Omega)}
					\to \tfrac{1}{\rho_{1}} \langle V^{3}, \phi\rangle_{L^{2}(\Omega)}. \notag
				\end{align}
				Therefore, $\tfrac{1}{\rho_{1}} V^{3} = F^{1}$,
				i.e., $(\mathcal{A} V)^{1} = F^{1}$.
				Taking into account $F^{1} \in H^{1}_{\Gamma_{1}}(\Omega)$, we conclude $V^{3} \in H^{1}_{\Gamma_{1}}(\Omega)$.
				
				\item[{\it \roman{zaehler})}] \addtocounter{zaehler}{1}
				Letting $\Phi = (0, \phi, 0, 0, 0, 0)$, $\phi \in (H^{1}_{\Gamma_{1}}(\Omega)$, we similarly get
				$(\mathcal{A} V)^{2} = F^{2}$ und $V^{4} \in (H^{1}_{\Gamma_{1}}(\Omega))^{2}$.
				
				\item[{\it \roman{zaehler})}] \addtocounter{zaehler}{1}
				Further, we choose $\Phi = (0, 0, \phi, 0, 0, 0)'$, $\phi \in H^{1}_{\Gamma_{1}}(\Omega)$.
				This yields
				\begin{align}
					\langle F^{3}, \phi\rangle_{L^{2}(\Omega)} \leftarrow &\tfrac{1}{\rho_{1}}
					\langle K \triangle V^{1}_{n} + K \mathrm{div}\, V^{2}_{n}, \phi\rangle_{L^{2}(\Omega)} \notag \\
					= &-\tfrac{K}{\rho_{1}} \langle \nabla V^{1}_{n}, \nabla \phi\rangle_{(L^{2}(\Omega))^{2}}
					+ \langle K \mathrm{div}\, V^{2}_{n}, \phi\rangle_{L^{2}(\Omega)} \notag \\
					\to &\tfrac{K}{\rho_{1}} \langle \nabla V^{1}, \nabla \phi\rangle_{(L^{2}(\Omega))^{2}}
					+ \langle K \mathrm{div}\, V^{2}, \phi\rangle_{L^{2}(\Omega)} \notag
				\end{align}
				implying $\triangle V^{1} \in L^{2}(\Omega)$ and
				$\tfrac{1}{\rho_{1}}(K \triangle V^{1} + K \mathrm{div}\, V^{2}) = F^{3}$, i.e., $(\mathcal{A} V)^{3} = F^{3}$.
				
				\item[{\it \roman{zaehler})}] \addtocounter{zaehler}{1}
				For $\Phi = (0, 0, 0, 0, 0, \phi)'$, $\phi \in (H^{1}_{\Gamma_{3}}(\Omega))^{2}$, we obtain
				\begin{align}
					\langle F^{6}, \phi\rangle_{(L^{2}(\Omega))^{2}} \leftarrow &\tfrac{1}{\tau_{0}}
					\langle -\kappa \nabla V^{5}_{n} - \delta V^{6}_{n}, \phi\rangle_{(L^{2}(\Omega))^{2}}
					= \tfrac{\kappa}{\tau_{0}} \langle V^{5}_{n}, \mathrm{div}\, \phi\rangle_{L^{2}(\Omega)}
					- \tfrac{\delta}{\tau_{0}} \langle V^{6}_{n}, \phi\rangle_{(L^{2}(\Omega))^{2}} \notag \\
					\to &\tfrac{\kappa}{\tau_{0}} \langle V^{5}, \mathrm{div}\, \phi\rangle_{L^{2}(\Omega)}
					- \tfrac{\delta}{\tau_{0}} \langle V^{6}, \phi\rangle_{(L^{2}(\Omega))^{2}}. \notag
				\end{align}
				Hence, $V^{5} \in H^{1}_{\Gamma_{3}}(\Omega)$ and
				$\tfrac{1}{\tau_{0}} (-\kappa \nabla V^{5} - \delta V^{6}) = F^{6}$,
				i.e., $(\mathcal{A} V)^{6} = F^{6}$.
				
				\item[{\it \roman{zaehler})}] \addtocounter{zaehler}{1}
				Selecting now $\Phi = (0, 0, 0, \phi, 0, 0)'$, $\phi \in (H^{1}_{\Gamma_{1}}(\Omega))^{2}$,
				we find
				\begin{align}
					\langle F^{4}, \phi\rangle_{(L^{2}(\Omega))^{2}} \leftarrow &\tfrac{1}{\rho_{2}}
					\langle -K \nabla V^{1}_{n} + \mathcal{D}' S \mathcal{D} V^{2}_{n} - K V^{2}_{n} - \gamma \nabla V^{5}_{n}, \phi\rangle_{(L^{2}(\Omega))^{2}} \notag \\
					= &\tfrac{K}{\rho_{2}} \langle S \mathcal{D} V^{2}_{n}, \mathcal{D} \phi\rangle_{(L^{2}(\Omega))^{3}}
					+ \tfrac{1}{\rho_{2}} \langle -K \nabla V^{1}_{n} - K V^{2}_{n} - \gamma \nabla V^{5}_{n}, \phi\rangle_{(L^{2}(\Omega))^{2}} \notag \\
					\to &\tfrac{K}{\rho_{2}} \langle S \mathcal{D} V^{2}, \mathcal{D} \phi\rangle_{(L^{2}(\Omega))^{3}}
					+ \tfrac{1}{\rho_{2}} \langle -K \nabla V^{1} - K V^{2} - \gamma \nabla V^{5}, \phi\rangle_{(L^{2}(\Omega))^{2}}. \notag
				\end{align}
				Thus, $\mathcal{D}' S \mathcal{D} V^{2} \in (L^{2}(\Omega))^{2}$ and
				$\tfrac{1}{\rho_{2}} (-K \nabla V^{1} + \mathcal{D}' S \mathcal{D} V^{2} - K V^{2} - \gamma \nabla V^{5}) = F^{4}$,
				i.e., $(\mathcal{A} V)^{4} = F^{4}$.
				
				\item[{\it \roman{zaehler})}] \addtocounter{zaehler}{1}
				Finally, we let $\Phi = (0, 0, 0, 0, \phi, 0)'$ mit $\phi \in H^{1}_{\Gamma_{3}}(\Omega)$ and deduce
				\begin{align}
					\langle F^{5}, \phi\rangle_{L^{2}(\Omega)} \leftarrow &\tfrac{1}{\rho_{3}}
					\langle -\gamma \mathrm{div}\, V^{4}_{n} - \beta V^{5}_{n} - \kappa \mathrm{div}\, V^{6}_{n}, \phi\rangle_{L^{2}(\Omega)} \notag \\
					= &\tfrac{\kappa}{\rho_{3}} \langle V^{6}_{n}, \nabla \phi\rangle_{L^{2}(\Omega)}
					- \tfrac{1}{\rho_{3}} \langle \gamma \mathrm{div}\, V^{4}_{n} + \beta V^{5}_{n}, \phi\rangle_{L^{2}(\Omega)} \notag \\
					\to &\tfrac{\kappa}{\rho_{3}} \langle V^{6}_{n}, \nabla \phi\rangle_{L^{2}(\Omega)}
					- \tfrac{1}{\rho_{3}} \langle \gamma \mathrm{div}\, V^{4} + \beta V^{5}, \phi\rangle_{L^{2}(\Omega)} \notag
				\end{align}
				implying that $\mathrm{div}\, V^{6} \in L^{2}(\Omega)$
				and $\tfrac{1}{\rho_{3}} (-\gamma \mathrm{div}\, V^{4} - \beta V^{5} - \kappa \mathrm{div}\, V^{6}) = F^{5}$ hold true,
				i.e., $(\mathcal{A} V)^{5} = F^{5}$.
			\end{itemize}
				
			There remains to show that $V$ satisfies the generalized Neumann boundary conditions (\ref{GLEICHUNG_NEUMANN_MR_LIN_KF_1})--(\ref{GLEICHUNG_NEUMANN_MR_LIN_KF_3}).
			To this end, we proceed as follows.
			\begin{itemize}
				\setcounter{zaehler}{1}
				\item[{\it \roman{zaehler})}] \addtocounter{zaehler}{1}
				Let $\phi \in H^{1}_{\Gamma_{1}}(\Omega)$. Then
				\begin{align}
					\langle \triangle V^{1} + \mathrm{div}\, V^{2}, \phi\rangle_{L^{2}(\Omega)}
					\leftarrow &\langle \triangle V^{1}_{n} + \mathrm{div}\, V^{2}_{n}, \phi\rangle_{L^{2}(\Omega)} \notag \\
					= &-\langle \nabla V^{1}_{n} + V^{2}_{n}, \nabla \phi \rangle_{(L^{2}(\Omega))^{2}}
					\to \langle \nabla V^{1} + V^{2}, \nabla \phi \rangle_{(L^{2}(\Omega))^{2}}. \notag
				\end{align}
			
				\item[{\it \roman{zaehler})}] \addtocounter{zaehler}{1}
				For $\phi \in (H^{1}_{\Gamma_{1}}(\Omega))^{2}$, we get
				\begin{align}
					\langle \mathcal{D}' S \mathcal{D} V^{2} - \gamma \nabla V^{5}, \phi\rangle_{(L^{2}(\Omega))^{2}}
					\leftarrow &\langle \mathcal{D}' S \mathcal{D} V^{2}_{n} - \gamma \nabla V^{5}_{n}, \phi\rangle_{(L^{2}(\Omega))^{2}} \notag \\
					= &-\langle S \mathcal{D} V^{2}_{n}, \mathcal{D} \phi \rangle_{(L^{2}(\Omega))^{2}} + \gamma \langle V^{5}_{n}, \mathrm{div}\, \phi \rangle_{L^{2}(\Omega)} \notag \\
					\to &-\langle S \mathcal{D} V^{2}, \mathcal{D} \phi \rangle_{(L^{2}(\Omega))^{2}} + \gamma \langle V^{5}, \mathrm{div}\, \phi \rangle_{L^{2}(\Omega)}. \notag
				\end{align}
			
				\item[{\it \roman{zaehler})}] \addtocounter{zaehler}{1}
				Choosing an arbitrary $\phi \in H^{1}_{\Gamma_{3}}(\Omega)$, we finally obtain
				\begin{align}
					\langle \mathrm{div}\, V^{6}, \phi\rangle_{L^{2}(\Omega)}
					\leftarrow & \langle \mathrm{div}\, V^{6}_{n}, \phi\rangle_{(L^{2}(\Omega))^{2}}
					= \langle V^{6}_{n}, \nabla \phi\rangle_{L^{2}(\Omega)}
					\to -\langle V^{6}, \nabla \phi\rangle_{(L^{2}(\Omega))^{2}}. \notag
				\end{align}
			\end{itemize}
			Alltogether, we have shown that $\mathcal{A}$ is a closed operator.
				
			\item
			Next, we show $\mathrm{im}(\lambda - \mathcal{A}) = \mathcal{H}$ for all $\lambda > 0$.
			To this end, we prove that the equation
			\begin{equation}
				(\lambda - \mathcal{A}) V = F \label{REISSNER_MINDLIN_ELLIPTISCHES_PARAMETERABHAENGIGES_PROBLEM}
			\end{equation}
			is solvable for any $F \in \mathcal{H}$.
			Since $D(\mathcal{A})$ is a dense subset of $\mathcal{H}$ and $\mathcal{A}$ is closed, we can select $F \in D(\mathcal{A})$.	
			Thus, for $F \in D(\mathcal{A})$, we are looking for solutions of
			\begin{equation}
				\begin{split}
					\lambda V^{1} - V^{3} &= F^{1}, \\
					\lambda V^{2} - V^{4} &= F^{2}, \\
					\lambda V^{3} - K \triangle V^{1} - K \mathrm{div}\, V^{2} &= \rho_{1} F^{3}, \\
					\lambda V^{4} + K \nabla V^{1} - \mathcal{D}' S \mathcal{D} V^{2} + K V^{2} + \gamma \nabla V^{5} &= \rho_{2} F^{4}, \\
					\lambda V^{5} + \gamma \mathrm{div}\, V^{4} + \beta V^{5} + \kappa \mathrm{div}\, V^{6} &= \rho_{3} F^{5}, \\
					\lambda V^{6} + \kappa \nabla V^{5} + \delta V^{6} &= \tau_{0} F^{6}.
				\end{split} \notag
			\end{equation}
			To eliminate $V^{3}$, $V^{4}$, we substitute
			\begin{align}
				V^{3} &= \lambda V^{1} - F^{1}, &
				V^{4} &= \lambda V^{2} - F^{2}, &
				V^{6} &= \tfrac{1}{\lambda + \delta} (- \kappa \nabla V^{5} + \tau_{0} F^{6}) \notag
			\end{align}
			and obtain
			\begin{equation}
				\begin{split}
					\lambda (\lambda + d) V^{1} - K \triangle V^{1} - K \mathrm{div}\, V^{2} &= G_{1}, \\
					\lambda^{2} V^{2} + K \nabla V^{1} - \mathcal{D}' S \mathcal{D} V^{2} + K V^{2} + \gamma \nabla V^{5} &= G_{2}, \\
					\lambda V^{5} + \gamma \lambda \mathrm{div}\, V^{2} + \beta V^{5} - \tfrac{\kappa^{2}}{1 + \delta} \triangle V^{5} &= G_{3}
				\end{split} 
				\label{REISSNER_MINDLIN_REDUZIERTES_ELIPTISCHES_PROBLEM}
			\end{equation}
			with
			\begin{align}
				G_{1} = \rho_{1} F^{3} + \lambda F^{1}, \quad
				G_{2} = \rho_{2} F^{4} + \lambda F^{2}, \quad
				G_{3} = \rho_{3} F^{5} + \gamma \mathrm{div}\, F^{2} + \tfrac{\tau_{0} \kappa}{\lambda + \delta} \mathrm{div}\, F^{6}. \notag
			\end{align}
		
			To solve the elliptic problem (\ref{REISSNER_MINDLIN_REDUZIERTES_ELIPTISCHES_PROBLEM}),
			we exploit the lemma of Lax \& Milgram.
			We consider the Hilbert space
			\begin{equation}
				\mathcal{V} := H^{1}_{\Gamma_{1}}(\Omega) \times (H^{1}_{\Gamma_{1}}(\Omega))^{2} \times H^{1}_{\Gamma_{3}}(\Omega) \notag
			\end{equation}
			equipped with the standard norm and introduce the bilinear form $a \colon \mathcal{V} \times \mathcal{V} \to \mathbb{R}$ via
			\begin{equation}
				\begin{split}
					a(V, W) := &
					\lambda^{3} \langle V^{1}, W^{1}\rangle_{L^{2}(\Omega)}
					+ \lambda^{3} \langle V^{2}, W^{2}\rangle_{(L^{2}(\Omega))^{2}}
					+ (\lambda + \beta) \langle V^{5}, W^{5}\rangle_{L^{2}(\Omega)} \\
					&K \lambda \langle \nabla V^{1} + V^{2}, \nabla W^{1} + W^{2}\rangle_{(L^{2}(\Omega))^{2}}
					+ \lambda \langle S \mathcal{D} V^{2}, \mathcal{D} W^{2}\rangle_{(L^{2}(\Omega))^{3}} + \\
					&\tfrac{\kappa^{2}}{\lambda + \delta} \langle \nabla V^{5}, \nabla W^{5}\rangle_{(L^{2}(\Omega))^{2}}
					+ \gamma \lambda \langle \nabla V^{5}, W^{2}\rangle_{(L^{2}(\Omega))^{2}}
					+ \gamma \lambda \langle \mathrm{div}\, V^{2}, W^{5}\rangle_{L^{2}(\Omega)}.
				\end{split} \label{REISSNER_MINDLIN_REDUZIERTES_ELIPTISCHES_PROBLEM_SCHWACH}
			\end{equation}
			After multiplying the equations in (\ref{REISSNER_MINDLIN_REDUZIERTES_ELIPTISCHES_PROBLEM})
			scalar in $L^{2}(\Omega)$, $(L^{2}(\Omega))^{2}$ and $L^{2}(\Omega)$ with $\lambda V^{1}$, $\lambda V^{2}$ and $V^{3}$, respectively,
			summing up the resulting equations and performing a partial integration,
			we obtain a weak formulation of Equation (\ref{REISSNER_MINDLIN_REDUZIERTES_ELIPTISCHES_PROBLEM}) in the form:
			Determine $V \in \mathcal{V}$ such that
			\begin{equation}
				a(V, W) = \lambda \langle G^{1}, W^{1}\rangle_{L^{2}(\Omega)} + \lambda \langle G^{2}, W^{2}\rangle_{(L^{2}(\Omega))^{2}}
				+ \langle G^{3}, W^{5}\rangle_{L^{2}(\Omega)} \notag
			\end{equation}
			for any $W \in \mathcal{V}$.

			The bilinear form $a$ is continuous and coercive on $\mathcal{V}$ due to the boundary conditions and
			the Korn's inequality from Theorem \ref{SATZ_KORNSCHE_UNGLEICHUNG}.
			The functional
			\begin{equation}
				\mathcal{V} \ni W \mapsto \lambda \langle G^{1}, W^{1}\rangle_{L^{2}(\Omega)} + \lambda \langle G^{2}, W^{2}\rangle_{(L^{2}(\Omega))^{2}}
				+ \langle G^{3}, W^{5}\rangle_{L^{2}(\Omega)} \notag
			\end{equation}
			is linear and continuous on $\mathcal{V}$.
			Applying now lemma of Lax \& Milgram, we deduce the existence of a weak solution $V \in \mathcal{V}$ to (\ref{REISSNER_MINDLIN_REDUZIERTES_ELIPTISCHES_PROBLEM_SCHWACH})
			which, in its turn, solves (\ref{REISSNER_MINDLIN_REDUZIERTES_ELIPTISCHES_PROBLEM}), too.
		
			Letting
			\begin{align}
				V^{3} &= \lambda V^{1} - F^{1}, &
				V^{4} &= \lambda V^{2} - F^{2}, &
				V^{6} &= \tfrac{1}{\lambda + \delta} (- \kappa \nabla V^{5} + \tau_{0} F^{6}), \notag
			\end{align}
			we conclude that $V = (V^{1}, \dots, V^{6})'$ solves Equation (\ref{REISSNER_MINDLIN_ELLIPTISCHES_PARAMETERABHAENGIGES_PROBLEM}).

			Thus, we have shown that $D(\mathcal{A}) \subset \mathrm{im}(\lambda - \mathcal{A})$.
			Since $D(\mathcal{A})$ is dense in $\mathcal{H}$ and $\mathrm{im}(\mathcal{A})$ is closed in $\mathcal{H}$, we finally obtain $\mathrm{im}(\lambda - \mathcal{A}) = \mathcal{H}$.
		\end{enumerate}
	\end{proof}
	
	We can now apply the theorem of Lumer \& Phillips to the Cauchy problem (\ref{REISSNER_MINDLIN_EVOLUTIONSPROBLEM}) to obtain the following existence result.
	\begin{theorem} \label{SATZ_WOHLGESTELLTHEIT}
		Let $V_{0} \in D(\mathcal{A})$.
		There exists then a unique classical solution to Equation (\ref{REISSNER_MINDLIN_EVOLUTIONSPROBLEM}) satisfying
		\begin{equation}
			V \in \mathcal{C}^{1}([0, \infty), \mathcal{H}) \cap \mathcal{C}^{0}([0, \infty), D(\mathcal{A})). \notag
		\end{equation}
		Moreover, if $V_{0} \in D(\mathcal{A}^{s})$ for a certain $s \in \mathbb{N}$, then we additionally have
		\begin{equation}
			V \in \bigcap_{k=0}^{s} \mathcal{C}^{k}([0, \infty), D(\mathcal{A}^{s-k})), \notag
		\end{equation}
		where $D(\mathcal{A}^{0}) := \mathcal{H}$.
	\end{theorem}

\section{Exponential stability}
	In this section, we study the stability properties of Equations (\ref{GLEICHUNG_REISNNER_MINDLIN_LINEAR_KOMPAKTE_FORM_1})--(\ref{GLEICHUNG_REISNNER_MINDLIN_LINEAR_KOMPAKTE_FORM_4})
	subject to Dirichlet boundary conditions for the elastic part and Neumann boundary conditions for the thermal part
	in two situations.
	First, we look at the case of a frictional damping on all elastic variables.
	Second, we restrict ourselves to the rotationally symmetric situation
	but retain only the frictional damping for the bending component $w$.

	For a number $d \geq 0$ and a symmetric, positive semidefinite matrix $D \in \mathbb{R}^{3 \times 3}$, we consider thus the problem
	\begin{align}
		\rho_{1} w_{tt} - K \mathrm{div}\, (v + \nabla w) + d w_{t} &= 0 \text{ in } (0, \infty) \times \Omega,
		\label{GLEICHUNG_REISNNER_MINDLIN_LINEAR_MECHANISCH_GEDAEMPFT_1} \\
		\rho_{2} v_{tt} - \mathcal{D}' S \mathcal{D} v + K (v + \nabla w) + \gamma \nabla \theta + D v_{t} &= 0 \text{ in } (0, \infty) \times \Omega,
		\label{GLEICHUNG_REISNNER_MINDLIN_LINEAR_MECHANISCH_GEDAEMPFT_2} \\
		\rho_{3} \theta_{t} + \kappa \mathrm{div}\, q + \beta \theta + \gamma \mathrm{div}\, v_{t} &= 0 \text{ in } (0, \infty) \times \Omega,
		\label{GLEICHUNG_REISNNER_MINDLIN_LINEAR_MECHANISCH_GEDAEMPFT_3} \\
		\tau_{0} q_{t} + \delta q + \kappa \nabla \theta &= 0 \text{ in } (0, \infty) \times \Omega
		\label{GLEICHUNG_REISNNER_MINDLIN_LINEAR_MECHANISCH_GEDAEMPFT_4}
	\end{align}
	subject to the boundary conditions
	\begin{align}
		w = |v| &= 0 \text{ on } (0, \infty) \times \Gamma, 
		\label{GLEICHUNG_REISNNER_MINDLIN_LINEAR_MECHANISCH_GEDAEMPFT_RB_1} \\
		q \cdot \nu &= 0 \text{ on } (0, \infty) \times \Gamma
		\label{GLEICHUNG_REISNNER_MINDLIN_LINEAR_MECHANISCH_GEDAEMPFT_RB_2}
	\end{align}
	and the initial conditions
	\begin{equation}
		w(0, \cdot) = w^{0}, \; w_{t}(0, \cdot) = w^{1}, \;
		v(0, \cdot) = v^{0}, \; v_{t}(0, \cdot) = v^{1}, \;
		\theta(0, \cdot) = \theta^{0}, \; q(0, \cdot) = q^{0},
		\label{GLEICHUNG_REISNNER_MINDLIN_LINEAR_MECHANISCH_GEDAEMPFT_AB}
	\end{equation}
	Despite of the notation abuse,
	the matrix $D \in \mathrm{SPD}(\mathbb{R}^{3})$ should not be confused with constant $D > 0$ consituting the matrix $S$.
	The natural first order energy associated with (\ref{GLEICHUNG_REISNNER_MINDLIN_LINEAR_MECHANISCH_GEDAEMPFT_1})--(\ref{GLEICHUNG_REISNNER_MINDLIN_LINEAR_MECHANISCH_GEDAEMPFT_4}) reads as
	\begin{equation}
		\mathcal{E}(t) := \tfrac{\rho_{1}}{2} \|w_{t}\|_{L^{2}(\Omega)}^{2} + \tfrac{\rho_{2}}{2} \|v_{t}\|_{(L^{2}(\Omega))^{2}}^{2}
		+ \tfrac{1}{2} \|\sqrt{S} \mathcal{D} v\|_{(L^{2}(\Omega))^{3}}^{2} + \tfrac{K}{2} \|v + \nabla w\|_{(L^{2}(\Omega))^{2}}^{2} +
		\tfrac{\rho_{3}}{2} \|\theta\|_{L^{2}(\Omega)}^{2} + \tfrac{\tau_{0}}{2} \|q\|_{(L^{2}(\Omega))^{2}}. \notag
	\end{equation}

\subsection{Full mechanical and thermal damping}
	First, we address the case of a full mechanical and thermal damping, i.e., $d > 0$, $D \in \mathrm{SPD}(\mathbb{R}^{2})$, $\beta > 0$.
	Analogous results for the equations of thermoelasticity with a mechanical damping were proved by Racke in \cite{Ra1993}
	for the case of parabolic heat conduction and by Ritter in \cite{Ri2011} for the case of hyperbolic heat conduction due to Cattaneo. 
	\begin{theorem}
		Let the parameters satisfy $\rho_{1}, \rho_{2}, \rho_{3}, \tau_{0}, K, \kappa, \delta, \gamma, d > 0, \beta > 0$,
		$S \in \mathrm{SPD}(\mathbb{R}^{3})$, $D \in \mathrm{SPD}(\mathbb{R}^{2})$.
		There exist then positive constants $C$ and $\alpha$ such that
		\begin{equation}
			\mathcal{E}(t) \leq C \mathcal{E}(0) e^{-2 \alpha t} \notag
		\end{equation}
		holds true for all $t \geq 0$.
		The latter depend neither on the initial data, nor on $t$
		and can be explicitely estimated based on the parameters and the domain $\Omega$.
	\end{theorem}

	\begin{proof}
		To prove the theorem, we want to construct a Lyapunov functional $\mathcal{F}$.
		Multiplying Equations (\ref{GLEICHUNG_REISNNER_MINDLIN_LINEAR_MECHANISCH_GEDAEMPFT_1}) and (\ref{GLEICHUNG_REISNNER_MINDLIN_LINEAR_MECHANISCH_GEDAEMPFT_3}) in $L^{2}(\Omega)$ with $w_{t}$ and $\theta$, respectively,
		as well as Equations (\ref{GLEICHUNG_REISNNER_MINDLIN_LINEAR_MECHANISCH_GEDAEMPFT_2}) and (\ref{GLEICHUNG_REISNNER_MINDLIN_LINEAR_MECHANISCH_GEDAEMPFT_4}) in $(L^{2}(\Omega))^{2}$ with $v_{t}$ and $q$, respectively,
		and exploiting the boundary conditions (\ref{GLEICHUNG_REISNNER_MINDLIN_LINEAR_MECHANISCH_GEDAEMPFT_RB_1}), (\ref{GLEICHUNG_REISNNER_MINDLIN_LINEAR_MECHANISCH_GEDAEMPFT_RB_2}),
		we find after a partial integration
		\begin{equation}
			\partial_{t} \mathcal{E}(t) \leq d \int_{\Omega} w_{t}^{2} \mathrm{d} x - \lambda \int_{\Omega} |v_{t}|^{2} \mathrm{d} x 
			-\beta \int_{\Omega} \theta^{2} \mathrm{d}x - \delta \int_{\Omega} |q|^{2} \mathrm{d} x
			\label{GLEICHUNG_REISSNER_MINDLIN_MECHANISCH_VOLLGEDAEMPFT_ENERGIE_ABGELEITET}
		\end{equation}
		with $\lambda := \min \sigma(D) > 0$ denoting the smallest eigenvalue of $D$.
		The function $\mathcal{F}$ has thus to be constructed in a way such that $\partial_{t} \mathcal{F}$ contains a negative multiple of $\mathcal{E}$,
		in particular, the terms $\int_{\Omega} |\nabla w|^{2} \mathrm{d} x$, $\int_{\Omega} |\sqrt{S} \mathcal{D} v|^{2} \mathrm{d} x$ and $\int_{\Omega} |\theta|^{2} \mathrm{d} x$.
		We define
		\begin{equation}
			\mathcal{F}_{1}(t) := \rho_{1} \int_{\Omega} w_{t} w \mathrm{d} x, \quad
			\mathcal{F}_{2}(t) := \rho_{1} \int_{\Omega} v_{t} \cdot v \mathrm{d} x \notag
		\end{equation}
		with $\cdot$ denoting the standard dot product on $\mathbb{R}^{2}$ and exploit Equations
		(\ref{GLEICHUNG_REISNNER_MINDLIN_LINEAR_MECHANISCH_GEDAEMPFT_1}), (\ref{GLEICHUNG_REISNNER_MINDLIN_LINEAR_MECHANISCH_GEDAEMPFT_2}) und (\ref{GLEICHUNG_REISNNER_MINDLIN_LINEAR_MECHANISCH_GEDAEMPFT_RB_1})
		to find after a partial integration
		\begin{equation}
			\begin{split}
				\partial_{t} \mathcal{F}_{1}(t) &=
				\int_{\Omega} (K \mathrm{div}\,\, (\nabla w + v) - d w_{t}) w \mathrm{d} x
				+ \rho_{1} \int_{\Omega} w_{t}^{2} \mathrm{d} x \\
				&=
				\int_{\Omega} - K (\nabla w + v) \cdot \nabla w \mathrm{d} x - d w_{t} w + \rho_{1} w_{t}^{2} \mathrm{d} x, \\
				\partial_{t} \mathcal{F}_{2}(t) &=
				\int_{\Omega} (\mathcal{D}'S\mathcal{D} v - K(v + \nabla w) - \gamma \nabla \theta_{t} - D v_{t}) \cdot v \mathrm{d} x
				+ \rho_{2} \int_{\Omega} |v_{t}|^{2} \mathrm{d} x \\
				&= \int_{\Omega} -|\sqrt{S} \mathcal{D} v|^{2} - K(v + \nabla w) \cdot v
				+ \gamma \theta_{t} \mathrm{div}\,\, v - D v_{t} \cdot v + \rho_{2} |v_{t}|^{2} \mathrm{d} x.
			\end{split}
			\label{GLEICHUNG_REISSNER_MINDLIN_LINEAR_MECHANISCH_VOLLGEDAEMPFT_HILFFUNKTIONALE_ABGELEITET}
		\end{equation}
		Using now Young's inequality, the first Poincar\'{e}'s and well as Korn's inequality,
		we can estimate for arbitrary $\varepsilon, \varepsilon' > 0$
		the functionals in (\ref{GLEICHUNG_REISSNER_MINDLIN_LINEAR_MECHANISCH_VOLLGEDAEMPFT_HILFFUNKTIONALE_ABGELEITET}) as follows:
		\begin{equation}
			\begin{split}
				\partial_{t} \mathcal{F}_{1}(t) &\leq
				\int_{\Omega} - K |\nabla w|^{2} + \tfrac{K}{2} |\nabla w|^{2} + \tfrac{K}{2} |v|^{2}
				+ \tfrac{d \varepsilon}{2} w^{2} + \big(\tfrac{d}{2\varepsilon} + \rho_{1}\big) w_{t}^{2} \mathrm{d} x \\
				&\leq
				\int_{\Omega} -\big(\tfrac{K}{2} - \tfrac{C_{\mathcal{P}} d \varepsilon}{2}) |\nabla w|^{2} + \tfrac{K}{2} |v|^{2}
				+ \big(\tfrac{d}{2\varepsilon} + \rho_{1}\big) w_{t}^{2} \mathrm{d} x, \\
				\partial_{t} \mathcal{F}_{2}(t) &\leq
				\int_{\Omega} -|\sqrt{S} \mathcal{D} v|^{2} - K |v|^{2} + \tfrac{K(1 + \varepsilon')}{2} |v|^{2} + \tfrac{K}{2(1 + \varepsilon')} |\nabla w|^{2}
				+ \tfrac{\gamma \varepsilon}{2} |\mathrm{div}\,\, v|^{2} \\
				&\phantom{\leq \int_{\Omega}} + \tfrac{\gamma}{2 \varepsilon} \theta^{2}
				+ \tfrac{\|D\| \varepsilon}{2} |v|^{2}
				+ \big(\tfrac{\|D\|}{2 \varepsilon} + \rho_{2}\big) |v_{t}|^{2} \mathrm{d} x \\
				&\leq
				\int_{\Omega}
				-\big(1 - \tfrac{K \varepsilon'}{2 C_{\mathcal{K}, 1}} - \tfrac{(\gamma + \|D\|) \varepsilon}{2 C_{\mathcal{K}, 1}}\big) |\sqrt{S} \mathcal{D} v|^{2}
				- \tfrac{K}{2} |v|^{2} + \tfrac{K}{2(1 + \varepsilon')} |\nabla w|^{2} \\
				&\phantom{\leq \int_{\Omega}} + \tfrac{\gamma}{2 \varepsilon} \theta^{2}
				+ \big(\tfrac{\|D\|}{2 \varepsilon} + \rho_{2}\big) |v_{t}|^{2} \mathrm{d} x, \\
			\end{split}
			\label{GLEICHUNG_REISSNER_MINDLIN_LINEAR_MECHANISCH_VOLLGEDAEMPFT_HILFFUNKTIONALE_ABGELEITET_UND_ABGESCHAETZT}
		\end{equation}
		where $C_{\mathcal{P}}$ denotes the Poincar\'{e}'s constant and $C_{\mathcal{K}, 1}$ stands for the Korn's constant from Lemma \ref{SATZ_KORNSCHE_UNGLEICHUNG}.
		We let
		\begin{equation}
			\mathcal{F}(t) := \mathcal{F}_{1}(t) + \mathcal{F}_{2}(t) + N \mathcal{E}(t) \notag
		\end{equation}
		and combine Equations
		(\ref{GLEICHUNG_REISSNER_MINDLIN_MECHANISCH_VOLLGEDAEMPFT_ENERGIE_ABGELEITET}) and (\ref{GLEICHUNG_REISSNER_MINDLIN_LINEAR_MECHANISCH_VOLLGEDAEMPFT_HILFFUNKTIONALE_ABGELEITET_UND_ABGESCHAETZT})
		to obtain
		\begin{equation}
			\begin{split}
				\partial_{t} \mathcal{F}(t)
				&\leq
				C_{w_{t}} \int_{\Omega} w_{t} \mathrm{d} x +
				C_{v_{t}} \int_{\Omega} |v_{t}|^{2} \mathrm{d} x +
				C_{\vartheta} \int_{\Omega} \vartheta_{t}^{2} \mathrm{d} x +
				C_{q} \int_{\Omega} |q|^{2} \mathrm{d} x +
				\\
				&\phantom{\leq\;\;}
				C_{\nabla w} \int_{\Omega} |\nabla w| \mathrm{d} x +
				C_{\sqrt{S} \mathcal{D} v} \int_{\Omega} |\sqrt{S} \mathcal{D} v|^{2} \mathrm{d} x,
			\end{split} \notag
        	\end{equation}
		where
		\begin{align}
			C_{w_{t}} &=
			N d - \tfrac{d}{2\varepsilon} + \rho_{1}, &
			C_{v_{t}} &=
			N \lambda - \tfrac{\|D\|}{2 \varepsilon} + \rho_{2} - \tfrac{\gamma}{2 \varepsilon}, \notag \\
			C_{\vartheta} &=
			N \beta - \tfrac{\gamma}{2 \varepsilon} + \rho_{3}, &
			C_{q} &= N \delta, \label{GLEICHUNG_REISSNER_MINDLIN_MECHANISCH_VOLLGEDAEMPFT_STABILITAETSKONSTANTEN} \\
			C_{\nabla w} &=
			\big[\tfrac{K}{2} - \tfrac{K}{2(1 + \varepsilon')}\big] - \tfrac{C_{\mathcal{P}} d \varepsilon}{2}, &
			C_{\sqrt{S} \mathcal{D} v} &= \big[1 - \tfrac{K \varepsilon'}{2 C_{\mathcal{K}, 1}}\big] - \tfrac{(\gamma + \|D\|) \varepsilon}{2 C_{\mathcal{K}, 1}}. \notag
		\end{align}
		Now, we select $\varepsilon' > 0$ to be sufficiently small
		such that the terms in the brackets from Equation (\ref{GLEICHUNG_REISSNER_MINDLIN_MECHANISCH_VOLLGEDAEMPFT_STABILITAETSKONSTANTEN}) become positive.
		Further, we fix a small $\varepsilon > 0$ to assure for
		$C_{\nabla w} > 0$ and $C_{\sqrt{S} \mathcal{D} v} > 0$.
		Finally, we pick a sufficiently large $N > 0$ such that all constants in  (\ref{GLEICHUNG_REISSNER_MINDLIN_MECHANISCH_VOLLGEDAEMPFT_STABILITAETSKONSTANTEN}) become positive.
		Thus,
		\begin{equation}
			C_{\min} := \min\{C_{w_{t}}, C_{v_{t}}, C_{\theta}, C_{q}, C_{\nabla w}, C_{\sqrt{S} \mathcal{D} v}\} > 0. \notag
		\end{equation}
		Using now the Korn's inequality from Lemma \ref{SATZ_KORNSCHE_UNGLEICHUNG}, we obtain
		\begin{equation}
			\partial_{t} \mathcal{F}(t)
			\leq -2C_{\min} \cdot \tfrac{\min\big\{1, C_{\mathcal{K}}\big\}}{
			\max\{1, \rho_{1}, \rho_{2}, \rho_{3}, \tau{0}\}} \mathcal{E}(t) =: - \tilde{C} \mathcal{E}(t).
		\end{equation}
		Taking into account
		\begin{equation}
			|(\mathcal{F}_{1} + \mathcal{F}_{2})(t)|
			\leq \tfrac{\max\big\{1, \rho_{1}, \rho_{2}\big\}}{\min\{1, C_{\mathcal{K}}\}} \mathcal{E}(t)
			=: \hat{C} \mathcal{E}(t) \notag
		\end{equation}
		we conclude
		\begin{equation}
			\beta_{1} \mathcal{E}(t) \leq \mathcal{L}(t) \leq \beta_{2} \mathcal{E}(t) \text{ for } t \geq 0 \notag
		\end{equation}
		with
		$\beta_{1} = N - \hat{C}$, $\beta_{2} = N + \hat{C}$.
		If neccessary, we increase $N$ to make $\beta_{1}$ positive.
		Gronwall's inequality now yields
		\begin{equation}
			\mathcal{E}(t) \leq \tfrac{1}{\beta_{1}} \mathcal{L}(t) \leq
			\tfrac{1}{\beta_{1}} \mathcal{E}(0) e^{-\tfrac{C}{\beta_{2}} t}
			=: C \mathcal{E}(0) e^{-2 \alpha t} \text{ for all } t \geq 0 \notag
		\end{equation}
		with $C, \alpha > 0$.
		This means that $E$ decays exponentially.
	\end{proof}

	\begin{remark}
		As a matter of fact, the constant $\beta$ must be positive in physical settings.
		Assuming
		\begin{equation}
			\int_{\Omega} \theta_{0} \mathrm{d}x = 0 \notag
		\end{equation}
		and using the functional $\mathcal{F}_{4}$ from the proof of Theorem \ref{SATZ_REISSNER_MINDLIN_LINEAR_EXPONENTIELLE_STABILITAET},
		our arguments can easily be carried over to the case $\beta = 0$.
		In contrast to Ritter's approach in \cite{Ri2011},
		no second order energy is required.
	\end{remark}

\subsection{Lack of strong stability in smooth domains}
	To justify the necessity of a frictional damping for both $w$ and $v$,
	we prove next that Equations (\ref{GLEICHUNG_REISNNER_MINDLIN_LINEAR_MECHANISCH_GEDAEMPFT_1})--(\ref{GLEICHUNG_REISNNER_MINDLIN_LINEAR_MECHANISCH_GEDAEMPFT_AB})
	even lack a strong stability for $d > 0$ and $D = 0$
	when considered in a bounded domain $\Omega$ with a smooth boundary $\Gamma$.
	In this case, the domain $\Omega$ contains a ray of geometrical optics perpendicularly reflected from $\Gamma$
	and one could theoretically perform constructions similar to those in \cite{Da1968} or \cite{LeZu1999} to prove a non-uniform decay rate even for a bigger class of domains.
	For simplicity, we restrict ourselves to the case of a smooth boundary allowing for the definition of Helmholz projection.
	We will namely show the imposibility of stabilizing the solenoidal part of $v$.
	
	To avoid a trivial null space, we impose for simplicity the following boundary conditions:
	\begin{align}
		w = |v| &= 0 \text{ on } (0, \infty) \times \Gamma, 
		\label{GLEICHUNG_REISNNER_MINDLIN_LINEAR_MECHANISCH_TEILWEISE_GEDAEMPFT_RB_1} \\
		\theta &= 0 \text{ on } (0, \infty) \times \Gamma.
		\label{GLEICHUNG_REISNNER_MINDLIN_LINEAR_MECHANISCH_TEILWEISE_GEDAEMPFT_RB_2}
	\end{align}
	It should though be pointed out that a similar result would also hold
	under any natural boundary conditions on $w$, $\theta$ and $q$
	provided Dirichlet boundary conditions are imposed on $v$ on the whole of $\Gamma$.

	\begin{theorem} \label{SATZ_FEHLER_DER_EXPONENTIELLEN_STABILITAET}
		Let the boundary $\Gamma$ be of class $\mathcal{C}^{2}$ and let $D = 0$.
		Problem (\ref{GLEICHUNG_REISNNER_MINDLIN_LINEAR_MECHANISCH_GEDAEMPFT_1})--(\ref{GLEICHUNG_REISNNER_MINDLIN_LINEAR_MECHANISCH_GEDAEMPFT_4}), (\ref{GLEICHUNG_REISNNER_MINDLIN_LINEAR_MECHANISCH_GEDAEMPFT_AB}),
		(\ref{GLEICHUNG_REISNNER_MINDLIN_LINEAR_MECHANISCH_TEILWEISE_GEDAEMPFT_RB_1}), (\ref{GLEICHUNG_REISNNER_MINDLIN_LINEAR_MECHANISCH_TEILWEISE_GEDAEMPFT_RB_2})
		is not strongly stable, in particular, not uniformly stable.
	\end{theorem}
	
	\begin{proof}
		Equations (\ref{GLEICHUNG_REISNNER_MINDLIN_LINEAR_MECHANISCH_GEDAEMPFT_1})--(\ref{GLEICHUNG_REISNNER_MINDLIN_LINEAR_MECHANISCH_GEDAEMPFT_4}), (\ref{GLEICHUNG_REISNNER_MINDLIN_LINEAR_MECHANISCH_GEDAEMPFT_AB}),
		(\ref{GLEICHUNG_REISNNER_MINDLIN_LINEAR_MECHANISCH_TEILWEISE_GEDAEMPFT_RB_1}), (\ref{GLEICHUNG_REISNNER_MINDLIN_LINEAR_MECHANISCH_TEILWEISE_GEDAEMPFT_RB_2})
		can be rewritten in the evolution form.
		Theorem \ref{SATZ_WOHLGESTELLTHEIT} yields then the existence of
		unique solution $V = (w, v, w_{t}, v_{t}, \theta, q)'$ given as an application of the strongly continuous semigroup of linear bounded operators to the initial data.
		
		Now, we want to select the initial data such that the solution component $v$ remains solenoidal, i.e., $\mathrm{div}\, v = \mathrm{div}\, v_{t} = 0$.
		Since $\Gamma$ is smooth, there exists the Helmholtz-projection (cf. \cite{So2001})
		\begin{equation}
			P \colon (L^{2}(\Omega))^{2} \to L^{2}_{\sigma}(\Omega) \notag
		\end{equation}
		into the Hilbert space
		\begin{equation}
			L^{2}_{\sigma}(\Omega) = \{u \in (L^{2}(\Omega))^{2} \,|\,
			\langle u, \nabla \varphi\rangle_{(L^{2}(\Omega))^{2}} = 0 \text{ for all }
			\varphi \in L^{1}_{\mathrm{loc}}(\Omega) \text{ such that }
			\nabla \varphi \in (L^{2}(\Omega))^{2}\}. \notag
		\end{equation}
		$P$ is an orthogonal operator and $L^{2}_{\sigma}(\Omega)$ is closed.
		
		Applying the operator $P$ to Equation (\ref{GLEICHUNG_REISNNER_MINDLIN_LINEAR_MECHANISCH_GEDAEMPFT_2}) and exploiting the representation
		\begin{equation}
			\mathcal{D}' S \mathcal{D} v = D \tfrac{1 - \mu}{2} \triangle v + D \tfrac{1 + \mu}{2} \nabla \mathrm{div}\, v, \notag
		\end{equation}
		we obtain an equation for $u := P v$
		\begin{equation}
			\begin{split}
				\rho_{2} u_{tt} - D \tfrac{1 - \mu}{2} P \triangle u + K u &= 0 \;\, \text{ in } (0, \infty) \times \Omega, \\
				u &= 0 \;\, \text{ in } (0, \infty) \times \Gamma, \\
				u(t, \cdot) = u_{0} := P v_{0}, \; u_{t}(t, \cdot) &= u_{1} := P v_{1} \text{ in } \Omega.
			\end{split} \label{VERALLGEMEINERTE_KLEIN_GORDON_GLEICHUNG}
		\end{equation}
		Equation (\ref{VERALLGEMEINERTE_KLEIN_GORDON_GLEICHUNG}) has a strong resemblance to the Klein-Gordon-Equation
		with an unbounded selfadjoint Dirichlet-Stokes-Operator $D P \triangle$.
		We define the operator
		\begin{equation}
			\mathcal{A} \colon D(\mathcal{A}) \subset L^{2}_{\sigma}(\Omega) \to L^{2}_{\sigma}(\Omega), \quad
			u \mapsto D \tfrac{1 - \mu}{2} P \triangle u + K u, \notag
		\end{equation}
		where
		\begin{equation}
			D(\mathcal{A}) = (H^{2}(\Omega) \cap H^{1}_{0}(\Omega))^{2} \cap L^{2}_{\sigma}(\Omega). \notag
		\end{equation}
		It is known (see, e.g., \cite{Ke2010}) that the spectrum $\sigma(D \tfrac{1 + \mu}{2} P \triangle)$ of $D P \triangle$ purely consists of a discrete point spectrum
		\begin{equation}
			\sigma(D \tfrac{1 + \mu}{2} P \triangle) = \sigma_{p}(D \tfrac{1 + \mu}{2} P \triangle) = \{\lambda_{k} \,|\, k \in \mathbb{N}\} \notag
		\end{equation}
		with the eigenvalues $\lambda_{k}$, $k \in \mathbb{N}$, of finite multiplicity satisfying
		$0 < \lambda_{1} \leq \lambda_{k} \leq \lambda_{k+1} \to \infty$ as $k \to \infty$.
		Hence, $\sigma(\mathcal{A}) = \sigma_{p}(\mathcal{A}) = \{\mu_{k} \,|\, k \in \mathbb{N}\}$ with $\mu_{k} = \lambda_{k} + K$ for $k \in \mathbb{N}$.
		
		Let $\nu^{\ast} \in \sigma(\mathcal{A})$ and let $u^{\ast} \in D(\mathcal{A})$ be the eigenfunction corresponding to $\nu^{\ast}$ with $\|u^{\ast}\|_{(L^{2}(\Omega))^{2}} = 1$.
		We set $u_{0} := u^{\ast}$, $u_{1} := 0$ and find that
		\begin{equation}
			u(t) := \cos\big(\sqrt{\tfrac{\nu^{\ast}}{\rho_{2}}} t\big) v^{\ast}, \quad t \in \mathbb{R}, \notag
		\end{equation}
		is a solution of (\ref{VERALLGEMEINERTE_KLEIN_GORDON_GLEICHUNG}).
		The energy associated with $u$ reads as
		\begin{equation}
			\begin{split}
				\mathcal{E}_{1}(t) &= \rho_{2} \|v_{t}\|_{(L^{2}(\Omega))^{2}} +
				D \|\nabla v\|_{(L^{2}(\Omega))^{2}}
				= \nu^{\ast} \cos^{2}\big(\sqrt{\tfrac{\nu^{\ast}}{\rho_{2}}} t\big) + D \nu^{\ast} \cos^{2}\big(\sqrt{\tfrac{\nu^{\ast}}{\rho_{2}}} t\big) \\
				&= \nu^{\ast} (1 + D) \cos^{2}\big(\sqrt{\tfrac{\nu^{\ast}}{\rho_{2}}} t\big) \nrightarrow 0 \text{ for } t \to \infty.
			\end{split} \notag
		\end{equation}

		Thus, $(w, v, \theta, q)' = (0, \cos\big(\sqrt{\tfrac{\nu^{\ast}}{\rho_{2}}} t\big) v^{\ast}, 0, 0)'$
		is a solution of the original problem
		(\ref{GLEICHUNG_REISNNER_MINDLIN_LINEAR_MECHANISCH_GEDAEMPFT_1})--(\ref{GLEICHUNG_REISNNER_MINDLIN_LINEAR_MECHANISCH_GEDAEMPFT_AB}),
		(\ref{GLEICHUNG_REISNNER_MINDLIN_LINEAR_MECHANISCH_TEILWEISE_GEDAEMPFT_RB_1}), (\ref{GLEICHUNG_REISNNER_MINDLIN_LINEAR_MECHANISCH_TEILWEISE_GEDAEMPFT_RB_2})
		for the initial conditions
		\begin{equation}
			w^{0} = w^{1} = 0, \quad
			v^{0} = u^{\ast}, v^{1} = 0, \quad
			\theta^{0} = 0, \quad q^{0} = 0, \notag
		\end{equation}
		which satisfies
		\begin{equation}
			\begin{split}
				\mathcal{E}(t) &= \rho_{2} \|v_{t}\|_{(L^{2}(\Omega))^{2}}
				+ D \|\nabla v_{1}\|_{(L^{2}(\Omega))^{2}}^{2} + D \|\nabla v_{2}\|_{(L^{2}(\Omega))^{2}}^{2}
				+ D \tfrac{1 + \mu}{2} \|\mathrm{div}\, v\|_{L^{2}(\Omega)}^{2} \\
				&= \rho_{2} \|v_{t}\|_{(L^{2}(\Omega))^{2}}
				+ D \|\nabla v_{1}\|_{(L^{2}(\Omega))^{2}}^{2} + D \|\nabla v_{2}\|_{(L^{2}(\Omega))^{2}}^{2}
				= \mathcal{E}_{1}(t) \nrightarrow 0 \text{ for } t \to \infty, \notag
			\end{split}
		\end{equation}
		where $v = (v_{1}, v_{2})'$.
		Hence, $\mathcal{E}$ does not decay.
	\end{proof}

\subsection{Exponential stability for rotationally symmetric plates}
	As we have seen before, a single frictional damping for $w$ is not sufficiently strong to stabilize a thermoelastic Reissner-Mindlin-Timoshenko plate
	(\ref{GLEICHUNG_REISNNER_MINDLIN_LINEAR_MECHANISCH_GEDAEMPFT_1})--(\ref{GLEICHUNG_REISNNER_MINDLIN_LINEAR_MECHANISCH_GEDAEMPFT_4}) for general data.
	Motivated by Racke's result in \cite{Ra2003}, we reduce the problem to the rotationally symmetric case
	making thus the vector field $v$ irrotational.
	Though arguments similar to those of Jiang and Racke in \cite[Theorem 4.2]{JiaRa2000} and Racke in \cite{Ra2003} made for
	the system of classical or hyperbolic thermoelasticity could be adopted in our case,
	we decided to propose our own approach incorporating the Bogovski\u i operator and, to some extent, being a generalization of the method
	applied by Messaoudi et al. in \cite{MePoSa2009} to a one-dimensional Timoshenko-beam.
	In addition to its technical novelty, a direct benefit of our approach
	lies in the fact that we only need to consider a first and not a second order energy.
	We would like to mention that the spectral approach of Grobbelaar (cf. \cite{Gro2012, Gro2013}) seems also to be applicable to our problem.
	At the same time, we do not require the assumption of simple connectedness on $\Omega$.

	We study Equations (\ref{GLEICHUNG_REISNNER_MINDLIN_LINEAR_MECHANISCH_GEDAEMPFT_1})--(\ref{GLEICHUNG_REISNNER_MINDLIN_LINEAR_MECHANISCH_GEDAEMPFT_4}) subject to the boundary conditions
	\begin{align}
		w = |v| &= 0 \text{ on } (0, \infty) \times \Gamma, 
		\label{GLEICHUNG_REISNNER_MINDLIN_LINEAR_MECHANISCH_GEDAEMPFT_RADIAL_RB_1} \\
		q \cdot \nu &= 0 \text{ on } (0, \infty) \times \Gamma
		\label{GLEICHUNG_REISNNER_MINDLIN_LINEAR_MECHANISCH_GEDAEMPFT_RADIAL_RB_2}
	\end{align}
	and the initial conditions
	\begin{equation}
		w(0, \cdot) = w^{0}, \; w_{t}(0, \cdot) = w^{1}, \;
		v(0, \cdot) = v^{0}, \; v_{t}(0, \cdot) = v^{1}, \;
		\theta(0, \cdot) = \theta^{0}, \; q(0, \cdot) = q^{0}.
		\label{GLEICHUNG_REISNNER_MINDLIN_LINEAR_MECHANISCH_GEDAEMPFT_RADIAL_AB}
	\end{equation}
	For the solution given in Theorem \ref{SATZ_WOHLGESTELLTHEIT},
	we assume the vanishing mean value for $\theta$
	\begin{equation}
		\int_{\Omega} \theta \mathrm{d} x = 0 \text{ in } (0, \infty)
		\label{GLEICHUNG_REISNNER_MINDLIN_LINEAR_MECHANISCH_GEDAEMPFT_MITTELWERTBEDINGUNG_THETA}
	\end{equation}
	as well as the vanishing rotation for $v$
	\begin{equation}
		\mathrm{rot} v = \partial_{2} v_{1} - \partial_{1} v_{2} = 0 \text{ in } (0, \infty).
		\label{GLEICHUNG_REISNNER_MINDLIN_LINEAR_MECHANISCH_GEDAEMPFT_ROTATIONSFREIHEIT}
	\end{equation}

	\begin{theorem}
		\label{SATZ_REISSNER_MINDLIN_LINEAR_EXPONENTIELLE_STABILITAET}
		Let $\Omega \subset \mathbb{R}^{2}$ be a rotationally symmetric bounded domain.
		Let the parameters satisfy $\rho_{1}, \rho_{2}, \rho_{3}, \tau_{0}, K, \kappa, \delta, \gamma, d > 0$, $\beta \geq 0$, $D = 0$
		and the matrix $S$ come from Equation (\ref{GLEICHUNG_STEIFIGKEITSMATRIX}).
		Further, let the data $w^{0}, w^{1}, v^{0}, v^{1}, \theta^{0}, q^{0}$ be radially symmetric in the sense of \cite[Definition 4.4]{JiaRa2000}
		and satisfy
		\begin{equation}
			\int_{\Omega} \theta_{0} \mathrm{d} x = 0. \notag
		\end{equation}
		There exist then positive constants $C$ and $\alpha$ such that for the energy $E$
		\begin{equation}
			\mathcal{E}(t) \leq C \mathcal{E}(0) e^{-2 \alpha t} \notag
		\end{equation}
		holds true for all $t \geq 0$.
		The latter depend neither on the initial data, nor on $t$
		and can be explicitely estimated based on the parameters and the domain $\Omega$.
	\end{theorem}

	\begin{proof}
		Theorem \ref{SATZ_WOHLGESTELLTHEIT} applied for the case $\Gamma_{2} = \Gamma_{3} = \emptyset$ yields the existence of a unique classical solution.
		After a straighforward modification, \cite[Lemma 4.6]{JiaRa2000} implies that the solution remains rotationally symmetric for all times $t \geq 0$.

		Without loss of generality, we may assume $\beta = 0$.
		Indeed, denoting with $\mathcal{E}_{\beta}$ the natural energy associated with the system subject to some fixed initial conditions for a given $\beta \geq 0$
		and assuming the existence of constants $C$ and $\alpha$ independent of the initial data such that
		\begin{equation}
			\mathcal{E}_{0}(t) \leq C \mathcal{E}_{0}(0) e^{-2 \alpha t} \text{ for } t \geq 0, \notag
		\end{equation}
		we take into account
		\begin{equation}
			\partial_{t} \mathcal{E}_{\beta}(t) \leq \partial_{t}\mathcal{E}_{0}(t) - \beta \int_{\Omega} \theta^{2} \mathrm{d} x \leq \partial_{t} \mathcal{E}_{0}(t) \text{ for } t \geq 0 \notag
		\end{equation}
		as well as $\mathcal{E}_{\beta}(0) = \mathcal{E}_{0}(0)$ to conclude
		\begin{equation}
			\mathcal{E}_{\beta}(t) \leq \mathcal{E}_{0}(t) \leq C \mathcal{E}_{0}(0) e^{-2 \alpha t} \text{ for } t \geq 0. \notag
		\end{equation}
		Thus, we let $\beta = 0$.
		As already mentioned, some of the following steps are motivated by the one-dimensional proof of Messaoudi et al. from \cite{MePoSa2009}.
		
		Multiplying (\ref{GLEICHUNG_REISNNER_MINDLIN_LINEAR_MECHANISCH_GEDAEMPFT_1}) and (\ref{GLEICHUNG_REISNNER_MINDLIN_LINEAR_MECHANISCH_GEDAEMPFT_3}) in $L^{2}(\Omega)$ with $w_{t}$ and $\theta$, respectively,
		as well as (\ref{GLEICHUNG_REISNNER_MINDLIN_LINEAR_MECHANISCH_GEDAEMPFT_2}) and (\ref{GLEICHUNG_REISNNER_MINDLIN_LINEAR_MECHANISCH_GEDAEMPFT_4}) in $(L^{2}(\Omega))^{2}$ with $v_{t}$ and $q$, respectively,
		and employing integration by parts, we find
		\begin{equation}
			\partial_{t} \mathcal{E}(t) = -d \int_{\Omega} w_{t}^{2} \mathrm{d} x - \delta \int_{\Omega} |q|^{2} \mathrm{d} x. \notag
		\end{equation}
		
		With the solution $u \in H^{1}_{0}(\Omega)$ to the Poisson equation
		\begin{equation}
			\begin{split}
				-\triangle u &= \mathrm{div}\, v \text{ in } \Omega,  \\
				u &= 0 \text{ auf } \Gamma,
			\end{split} \notag
		\end{equation}
		we obtain
		\begin{equation}
			\int_{\Omega} |\nabla u|^{2} \mathrm{d} x = -\int_{\Omega} v \cdot \nabla u \mathrm{d} x. \notag
		\end{equation}
		Young's inequality further yields
		\begin{equation}
			\int_{\Omega} |\nabla u|^{2} \mathrm{d} x \leq \tfrac{1}{2} \int_{\Omega} |v|^{2} \mathrm{d} x + \tfrac{1}{2} \int_{\Omega} |\nabla u|^{2} \mathrm{d} x \notag
		\end{equation}
		and, therefore,
		\begin{equation}
			\int_{\Omega} |\nabla u|^{2} \mathrm{d} x \leq \int_{\Omega} |v|^{2} \mathrm{d} x. \label{GLEICHUNG_GRADIENTEN_ABSCHAETZUNG_GEGEN_L2_NORM_V}
		\end{equation}
		Similarly,
		\begin{equation}
			\int_{\Omega} |\nabla u_{t}|^{2} \mathrm{d} x \leq \int_{\Omega} |v_{t}|^{2} \mathrm{d} x. \label{GLEICHUNG_GRADIENTEN_ABSCHAETZUNG_GEGEN_L2_NORM_V_T}
		\end{equation}
		We define the functional
		\begin{equation}
			\mathcal{F}_{1}(t) :=
			\int_{\Omega} \left(\rho_{1} w_{t} u + \rho_{2} v_{t} v - \tfrac{\gamma \tau_{0}}{\kappa} v q\right) \mathrm{d} x. \notag
		\end{equation}
		Taking into account Equation (\ref{GLEICHUNG_REISNNER_MINDLIN_LINEAR_MECHANISCH_GEDAEMPFT_1}), we use partial integration to obtain
		\begin{align}
			\partial_{t} &\int_{\Omega} \rho_{1} w_{t} u \mathrm{d} x
			= \rho_{1} \int_{\Omega} (w_{tt} u + w_{t} u_{t}) \mathrm{d} x \notag \\
			&= \phantom{-} K \int_{\Omega} \mathrm{div}\,\, (\nabla w + v) \cdot u \mathrm{d} x
			- d \int_{\Omega} w_{t} u \mathrm{d} x + \rho_{1} \int_{\Omega} w_{t} u_{t} \mathrm{d} x \notag \\
			&= -K \int_{\Omega} (\nabla w + v) \cdot \nabla u \mathrm{d} x
			- d \int_{\Omega} w_{t} u \mathrm{d} x + \rho_{1} \int_{\Omega} w_{t} u_{t} \mathrm{d} x \notag \\
			&= \phantom{-} K \int_{\Omega} w \triangle u \mathrm{d} x - K \int_{\Omega} v \cdot \nabla u \mathrm{d} x
			- d \int_{\Omega} w_{t} u \mathrm{d} x + \rho_{1} \int_{\Omega} w_{t} u_{t} \mathrm{d} x \notag \\
			&= -K \int_{\Omega} w \mathrm{div}\,\, v \mathrm{d} x + K \int_{\Omega} |\nabla u|^{2} \mathrm{d} x
			- d \int_{\Omega} w_{t} u \mathrm{d} x + \rho_{1} \int_{\Omega} w_{t} u_{t} \mathrm{d} x. \notag
		\end{align}
		By the virtue of Equation (\ref{GLEICHUNG_REISNNER_MINDLIN_LINEAR_MECHANISCH_GEDAEMPFT_2}), we similarly get
		\begin{align}
			\partial_{t} &\int_{\Omega} \rho_{2} v_{t} \cdot v \mathrm{d} x
			= \rho_{2} \int_{\Omega} v_{tt} \cdot v \mathrm{d} x + \rho_{2} \int_{\Omega} |v_{t}|^{2} \mathrm{d} x \notag \\
			&= \phantom{-} \int_{\Omega} \mathcal{D}' S \mathcal{D} v \cdot v \mathrm{d} x - K \int_{\Omega} (v + \nabla w) \cdot v \mathrm{d} x
			- \gamma \int_{\Omega} \nabla \theta \cdot v \mathrm{d} x + \rho_{2} \int_{\Omega} |v_{t}|^{2} \mathrm{d} x \notag \\
			&= -\int_{\Omega} |\sqrt{S} \mathcal{D} v|^{2} \mathrm{d} x - K \int_{\Omega} |v|^{2} \mathrm{d} x + K \int_{\Omega} w \mathrm{div}\,\, v \mathrm{d} x
			- \gamma \int_{\Omega} \nabla \theta \cdot v \mathrm{d} x + \rho_{2} \int_{\Omega} |v_{t}|^{2} \mathrm{d} x \notag
		\end{align}
		as well as
		\begin{align}
			\partial_{t} &\int_{\Omega} -\tfrac{\gamma \tau_{0}}{\kappa} v \cdot q \mathrm{d} x
			= -\tfrac{\gamma \tau_{0}}{\kappa} \int_{\Omega} v_{t} \cdot q \mathrm{d} x
			+ \tfrac{\gamma \delta }{\kappa} \int_{\Omega} v \cdot q \mathrm{d} x + \gamma \int_{\Omega} v \cdot \nabla \theta \mathrm{d} x. \notag
		\end{align}
		Finally, we conclude
		\begin{align}
			\partial_{t} \mathcal{F}_{1}(t) = &
			K \int_{\Omega} |\nabla u|^{2} \mathrm{d} x - K \int_{\Omega} |v|^{2} \mathrm{d} x
			- d \int_{\Omega} w_{t} u \mathrm{d} x + \rho_{1} \int_{\Omega} w_{t} u_{t} \mathrm{d} x + \notag \\
			& \rho_{2} \int_{\Omega} |v_{t}|^{2} \mathrm{d} x - \int_{\Omega} |\sqrt{S} \mathcal{D} v|^{2} \mathrm{d} x -\tfrac{\gamma \tau_{0}}{\kappa} \int_{\Omega} v_{t} \cdot q \mathrm{d} x + \tfrac{\gamma \delta }{\kappa} \int_{\Omega} v \cdot q \mathrm{d} x. \notag
		\end{align}
		Using now the first Poincar\'{e}'s inequality, Young's inequality and Korn's inequality from Lemma \ref{SATZ_KORNSCHE_UNGLEICHUNG}
		as well as the estimates from Equations (\ref{GLEICHUNG_GRADIENTEN_ABSCHAETZUNG_GEGEN_L2_NORM_V}) and (\ref{GLEICHUNG_GRADIENTEN_ABSCHAETZUNG_GEGEN_L2_NORM_V_T}), we obtain
		\begin{align}
			\begin{split}
				\partial_{t} \mathcal{F}_{1}(t) \leq & \phantom{-} \tfrac{d}{2} \int_{\Omega} \left(\varepsilon_{1} u^{2} + \tfrac{1}{\varepsilon_{1}} w_{t}^{2}\right) \mathrm{d} x
				+ \tfrac{\rho_{1}}{2} \int_{\Omega} \left(\varepsilon_{1} u_{t}^{2} + \tfrac{1}{\varepsilon_{1}} w_{t}^{2}\right) \mathrm{d} x
				+ \rho_{2} \int_{\Omega} |v_{t}|^{2} \mathrm{d} x \\
				&- \int_{\Omega} |\sqrt{S} \mathcal{D} v|^{2} \mathrm{d} x
				+ \tfrac{\gamma \tau_{0}}{2 \kappa} \int_{\Omega} \left(\varepsilon_{1} |v_{t}|^{2} + \tfrac{1}{\varepsilon_{1}} |q|^{2}\right) \mathrm{d} x
				+ \tfrac{\gamma \delta}{2 \kappa} \int_{\Omega} \left(\varepsilon_{1} |v|^{2} + \tfrac{1}{\varepsilon_{1}} |q|^{2}\right) \mathrm{d} x \\
				\leq & \phantom{-} \tfrac{d}{2} \int_{\Omega} \left(\varepsilon_{1} \tfrac{C_{\mathcal{P}}}{C_{\mathcal{K}, 1}} |\sqrt{S} \mathcal{D} v|^{2}
				+ \tfrac{1}{\varepsilon_{1}} w_{t}^{2}\right) \mathrm{d} x
				+ \tfrac{\rho_{1}}{2} \int_{\Omega} \left(C_{\mathcal{P}} \varepsilon_{1} |v_{t}|^{2} + \tfrac{1}{\varepsilon} w_{t}^{2}\right) \mathrm{d} x
				+ \rho_{2} \int_{\Omega} |v_{t}|^{2} \mathrm{d} x \\
				&- \int_{\Omega} |\sqrt{S} \mathcal{D} v|^{2} \mathrm{d} x
				+ \tfrac{\gamma \tau_{0}}{2 \kappa} \int_{\Omega} \left(\varepsilon_{1} |v_{t}|^{2} + \tfrac{1}{\varepsilon_{1}} |q|^{2}\right) \mathrm{d} x
				+ \tfrac{\gamma \delta}{2 \kappa} \int_{\Omega} \left(\tfrac{\varepsilon_{1}}{C_{\mathcal{K}, 1}} |\sqrt{S} \mathcal{D} v|^{2}
				+ \tfrac{1}{\varepsilon_{1}} |q|^{2}\right) \mathrm{d} x \\
				\leq & \phantom{-} \tfrac{\rho_{1} + d}{2 \varepsilon_{1}} \int_{\Omega} w_{t}^{2} \mathrm{d} x
				+ \left[\rho_{2} - \tfrac{\varepsilon_{1}}{2} \left(\rho_{1} C_{\mathcal{P}}
				+ \tfrac{\gamma \tau_{0}}{\kappa}\right) \right] \int_{\Omega} |v_{t}|^{2} \mathrm{d} x \\
				& - \left[1 - \tfrac{\varepsilon_{1}}{2 C_{\mathcal{K}, 1}} \left(d C_{\mathcal{P}} + \tfrac{\gamma \delta}{\kappa} \right) \right] \int_{\Omega} |\sqrt{S} \mathcal{D} v|^{2} \mathrm{d} x
				+ \tfrac{\gamma (\tau_{0} + \delta)}{2 \kappa \varepsilon_{1}} \int_{\Omega} |q|^{2} \mathrm{d} x
			\end{split}
			\label{EQUATION_ABSCHAETZUNG_F1}
		\end{align}
		with the Poincar\'{e}'s constant $C_{\mathcal{P}} = C_{\mathcal{P}}(\Omega) > 0$ and an arbitrary small number $\varepsilon_{1} > 0$ to be fixed later.
		Here, we estimated
		\begin{equation}
			\int_{\Omega} |u|^{2} \mathrm{d} x \leq C_{\mathcal{P}} \int_{\Omega} |\nabla u|^{2} \mathrm{d} x
			\leq C_{\mathcal{P}} \int_{\Omega} |v|^{2} \mathrm{d} x
			\leq \tfrac{C_{\mathcal{P}}}{C_{\mathcal{K}, 1}} \int_{\Omega} |\sqrt{S} \mathcal{D} v|^{2} \mathrm{d} x. \notag
		\end{equation}
		
		Next, we consider the functional
		\begin{equation}
			\mathcal{F}_{2}(t) := \rho_{1} \int_{\Omega} w_{t} w \mathrm{d} x \notag
		\end{equation}
		and use Equation (\ref{GLEICHUNG_REISNNER_MINDLIN_LINEAR_MECHANISCH_GEDAEMPFT_1}) to find
		\begin{align}
			\partial_{t} \mathcal{F}_{2}(t) &= \rho_{1} \int_{\Omega} w_{t}^{2} \mathrm{d} x
			+ K \int_{\Omega} \mathrm{div}\,\, (v + \nabla w) \cdot w \mathrm{d} x - d \int_{\Omega} w_{t} w \mathrm{d} x \notag \\
			&= \rho_{1} \int_{\Omega} w_{t}^{2} \mathrm{d} x - K \int_{\Omega} |\nabla w|^{2} \mathrm{d} x - K \int_{\Omega} v \cdot \nabla w \mathrm{d} x
			- d \int_{\Omega} w_{t} w \mathrm{d} x. \notag
		\end{align}
		The latter can be estimated as
		\begin{align}
			\begin{split}
				\partial_{t} \mathcal{F}_{2}(t) &\leq -K \int_{\Omega} |\nabla w|^{2} \mathrm{d} x
				+ \tfrac{K}{2} \int_{\Omega} \left(\varepsilon_{2} |\nabla w|^{2} + \tfrac{1}{\varepsilon_{2}} |v|^{2} \right) \mathrm{d} x + \\
				&\phantom{= -} \tfrac{d}{2} \int_{\Omega} \left(\varepsilon_{2} w^{2} + \tfrac{1}{\varepsilon_{2}} w_{t}^{2}\right) \mathrm{d} x
				+ \rho_{1} \int_{\Omega} w_{t}^{2} \mathrm{d} x \\
				&\leq -\left(K - \tfrac{\varepsilon_{2} C_{\mathcal{P}}}{2} (K + d)\right) \int_{\Omega} |\nabla w|^{2} \mathrm{d} x
				+ \tfrac{K}{\varepsilon_{2} \mathcal{C}_{\mathcal{K}, 1}} \int_{\Omega} |\sqrt{S} \mathcal{D} v|^{2} \mathrm{d} x
				+ \left(\tfrac{d}{2 \varepsilon_{2}} + \rho_{1}\right) \int_{\Omega} w_{t}^{2} \mathrm{d} x.
			\end{split}
			\label{EQUATION_ABSCHAETZUNG_F2}
		\end{align}
		
		Exploiting the fact
		\begin{align}
			0 &= \rho_{3} \int_{\Omega} \theta_{t} \mathrm{d} x + \kappa \int_{\Omega} \mathrm{div}\,\, q \mathrm{d} x + \gamma \int_{\Omega} \mathrm{div}\,\, v_{t} \mathrm{d} x \notag \\
			&= \rho_{3} \partial_{t} \int_{\Omega} \theta \mathrm{d} x + \kappa \int_{\Gamma} q \cdot \nu \mathrm{d} x
			+ \gamma \int_{\Gamma} v_{t} \cdot \nu \mathrm{d} x \,=\, \rho_{3} \partial_{t} \int_{\Omega} \theta \mathrm{d} x, \notag
		\end{align}
		we easily see
		\begin{equation}
			\int_{\Omega} \theta(t, x) \mathrm{d} x \equiv \int_{\Omega} \theta(0, x) \mathrm{d} x = \int_{\Omega} \theta_{0}(x) \mathrm{d} x = 0. \notag
		\end{equation}
		This enables us to apply the second Poincar\'{e}'s inequality to $\theta$.
		Using now the definition of Bogowski\u i operator $\mathcal{B}_{\mathrm{rot}}$ from Theorem \ref{SATZ_BOGOWSKII_STETIGKEIT}, we introduce the following functional
		\begin{equation}
			\mathcal{F}_{3}(t) := \rho_{2} \rho_{3} \int_{\Omega} \mathcal{B}_{\mathrm{rot}} \theta \cdot v_{t} \mathrm{d} x. \notag
		\end{equation}
		Exploiting Equations (\ref{GLEICHUNG_REISNNER_MINDLIN_LINEAR_MECHANISCH_GEDAEMPFT_2}) and (\ref{GLEICHUNG_REISNNER_MINDLIN_LINEAR_MECHANISCH_GEDAEMPFT_3}), we obtain
		\begin{align}
			\partial_{t} \mathcal{F}_{3}(t) &= \rho_{2} \rho_{3} \int_{\Omega} \mathcal{B}_{\mathrm{rot}} \theta_{t} \cdot v_{t} \mathrm{d} x
			+ \rho_{2} \rho_{3} \int_{\Omega} \mathcal{B}_{\mathrm{rot}} \theta \cdot v_{tt} \mathrm{d} x \notag \\
			&= \rho_{2} \int_{\Omega} \mathcal{B}_{\mathrm{rot}} (-\kappa \mathrm{div}\,\, q - \gamma \mathrm{div}\,\, v_{t}) \cdot v_{t} \mathrm{d} x +
			\rho_{3} \int_{\Omega} \mathcal{B}_{\mathrm{rot}} \theta \cdot (\mathcal{D}' S \mathcal{D} v - K(v + \nabla w) - \gamma \nabla \theta) \mathrm{d} x \notag \\
			&= -\rho_{2} \kappa \int_{\Omega} (\mathcal{B}_{\mathrm{rot}} \mathrm{div}\,\, q) \cdot v_{t} \mathrm{d} x - \rho_{2} \gamma \int_{\Omega} |v_{t}|^{2} \mathrm{d} x
			- \rho_{3} \int_{\Omega} (\mathcal{D} \mathcal{B}_{\mathrm{rot}} \theta) \cdot (S \mathcal{D} v) \mathrm{d} x \notag \\
			&\phantom{=} -\rho_{3} K \int_{\Omega} \mathcal{B}_{\mathrm{rot}} \theta \cdot (v + \nabla w) \mathrm{d} x
			+ \rho_{3} \gamma \int_{\Omega} \mathrm{div}\,\, \mathcal{B}_{\mathrm{rot}} \theta \cdot \theta \mathrm{d} x \notag \\
			&= -\rho_{2} \kappa \int_{\Omega} (\mathcal{B}_{\mathrm{rot}} \mathrm{div}\,\, q) \cdot v_{t} \mathrm{d} x - \rho_{2} \gamma \int_{\Omega} |v_{t}|^{2} \mathrm{d} x
			- \rho_{3} \int_{\Omega} (\mathcal{D} \mathcal{B}_{\mathrm{rot}} \theta) \cdot (S \mathcal{D} v) \mathrm{d} x \notag \\
			&\phantom{=} -\rho_{3} K \int_{\Omega} \mathcal{B}_{\mathrm{rot}} \theta \cdot (v + \nabla w) \mathrm{d} x
			- \rho_{3} \gamma \int_{\Omega} \theta^{2} \mathrm{d} x. \notag
		\end{align}
		We would like to stress
		that the injectivity of Bogowski\u i operator was essential here
		for us to be able to reconstruct $v_{t}$ from $\mathcal{B}_{\mathrm{rot}} \mathrm{div}\,\, v_{t}$.
		In general, this is not possible unless the vector field is irrotational and vanishes on $\Gamma$.
		Using the Young's inequality and exploiting the continuity of Bogowski\u i operator, we can estimate
		\begin{equation}
			\begin{split}
				\partial_{t} \mathcal{F}_{3}(t) &\leq - \rho_{2} \gamma \int_{\Omega} |v_{t}|^{2} \mathrm{d} x
				+ \tfrac{\rho_{2} \kappa}{2} \int_{\Omega} \left(\varepsilon_{3} |v_{t}|^{2} + \tfrac{C'_{\mathcal{B}_{\mathrm{rot}}}}{\varepsilon_{3}} |q|^{2}\right) \mathrm{d} x \\
				&\phantom{=}  + \tfrac{\rho_{3} \|S\|}{2} \int_{\Omega} \left(\varepsilon'_{3} |S \mathcal{D} v|^{2} + \tfrac{C_{\mathcal{B}}}{\varepsilon'_{3}} \theta^{2}\right) \mathrm{d} x
				+ \tfrac{\rho_{3} K}{2} \int_{\Omega} \left(\varepsilon'_{3} |v|^{2} + \tfrac{C_{\mathcal{K}, 2} C_{\mathcal{B}}}{\varepsilon'_{3}} \theta^{2} \right) \mathrm{d} x \\
				&\phantom{=}  + \tfrac{\rho_{3} K}{2} \int_{\Omega} \left(\varepsilon'_{3} |\nabla w|^{2} + \tfrac{c_{\mathcal{B}}}{\varepsilon'_{3}} \theta^{2} \right) \mathrm{d} x
				+ \rho_{3} \gamma \int_{\Omega} \theta^{2} \mathrm{d} x \\
				&\leq -\left(\rho_{2} \gamma - \tfrac{\rho_{2} \kappa \varepsilon_{3}}{2}\right) \int_{\Omega} |v_{t}|^{2} \mathrm{d} x
				+ \left(\tfrac{\rho_{3} 2 \|S\| \varepsilon_{3}'}{2} + \tfrac{\rho_{3} K C_{P} \varepsilon_{3}'}{2 C_{\mathcal{K}, 1}} \right) \int_{\Omega} |\sqrt{S} \mathcal{D} v|^{2} \mathrm{d} x
				+ \tfrac{\rho_{3} K \varepsilon'_{3}}{2} \int_{\Omega} |\nabla w|^{2} \mathrm{d} x
			\end{split}
			\notag
		\end{equation}
		\begin{equation}
		 	\begin{split}
				&\phantom{=} \left(\rho_{3} \gamma - \tfrac{\rho_{3} (\|S\| + K) C_{\mathcal{B}}}{2 \varepsilon'_{3}}
				- \tfrac{\rho_{3} K C_{\mathcal{B}}}{2 \varepsilon'_{3}} \theta^{2}\right) \int_{\Omega} \theta^{2} \mathrm{d} x
				+ \tfrac{\rho_{2} \kappa C'_{\mathcal{B}_{\mathrm{rot}}}}{2 \varepsilon_{3}} \int_{\Omega} |q|^{2} \mathrm{d} x
			\end{split}
			\label{EQUATION_ABSCHAETZUNG_F3}
		\end{equation}
		for arbitrary positive $\varepsilon_{3}$ and $\varepsilon'_{3}$.
		The constants $C_{\mathcal{B}_{\mathrm{rot}}}$ and $C'_{\mathcal{B}_{\mathrm{rot}}}$ occuring above
		come from Theorem \ref{SATZ_BOGOWSKII_STETIGKEIT} and Theorem \ref{SATZ_BOGOWSKII_ABSCHAETZUNG}.
		
		Finally, we define
		\begin{equation}
			\mathcal{F}_{4}(t) := -\tau_{0} \rho_{3} \int_{\Omega} q \cdot \mathcal{B}_{\mathrm{rot}} \theta \mathrm{d} x \notag
		\end{equation}
		and obtain
		\begin{align}
			\partial_{t} \mathcal{F}_{4}(t) &= -\rho_{3} \int_{\Omega} (-\delta q - \kappa \nabla \theta) \cdot \mathcal{B}_{\mathrm{rot}} \theta \mathrm{d} x
			-\tau_{0} \int_{\Omega} q \cdot \mathcal{B}_{\mathrm{rot}} (-\kappa \mathrm{div}\,\, q - \gamma \mathrm{div}\,\, v_{t}) \mathrm{d} x \notag \\
			&= \phantom{-}\rho_{3} \delta \int_{\Omega} q \cdot \mathcal{B} \theta \mathrm{d} x - \rho_{3} \kappa \int_{\Omega} \theta^{2} \mathrm{d} x
			+ \tau_{0} \kappa \int_{\Omega} q \cdot \mathcal{B} \mathrm{div}\,\, q \mathrm{d} x
			+ \tau_{0} \gamma \int_{\Omega} q \cdot v_{t} \mathrm{d} x \notag
		\end{align}
		since
		\begin{equation}
			\int_{\Omega} \nabla \theta \cdot \mathcal{B}_{\mathrm{rot}} \theta \mathrm{d} x
			= -\int_{\Omega} \theta \mathrm{div}\,\, \mathcal{B}_{\mathrm{rot}} \theta \mathrm{d} x
			= -\int_{\Omega} \theta^{2} \mathrm{d} x. \notag
		\end{equation}
		This yields the estimate
		\begin{align}
			\begin{split}
				\partial_{t} \mathcal{F}_{4} &\leq - \rho_{3} \kappa \int_{\Omega} \theta^{2} \mathrm{d} x
				+ \tfrac{\rho_{3} \delta}{2} \int_{\Omega} \left(\varepsilon_{4} C_{\mathcal{B}_{\mathrm{rot}}} \theta^{2} + \tfrac{1}{\varepsilon_{4}} |q|^{2} \right) \mathrm{d} x \\
				&\phantom{=} + \tau_{0} \kappa (1 + C'_{\mathcal{B}_{\mathrm{rot}}}) \int_{\Omega} |q|^{2} \mathrm{d} x
				+ \tfrac{\tau_{0} \gamma}{2} \int_{\Omega} \left(\varepsilon'_{4} |v_{t}|^{2} + \tfrac{1}{\varepsilon'_{4}} |q|^{2}\right) \mathrm{d} x \\
				&= \left(-\rho_{3} \kappa + \tfrac{\varepsilon_{4} \rho_{3} \delta C_{\mathcal{B}_{\mathrm{rot}}}}{2} \right) \int_{\Omega} \theta^{2} \mathrm{d} x
				+ \tfrac{\varepsilon'_{4} \tau_{0} \gamma}{2} \int_{\Omega} |v_{t}|^{2} \mathrm{d} x \\
				&\phantom{=} + \left((1 + C'_{\mathcal{B}_{\mathrm{rot}}}) \tau_{0} \kappa + \tfrac{\rho_{3} \delta}{2 \varepsilon_{4}} + \tfrac{\tau_{0} \gamma}{2 \varepsilon'_{4}}\right)
				\int_{\Omega} |q|^{2} \mathrm{d} x.
			\end{split}
			\label{EQUATION_ABSCHAETZUNG_F4}
		\end{align}

		For positive $N, N_{4}$, we define the auxiliary functional $\mathcal{F}$ by the means of
		\begin{equation}
			\mathcal{F}(t) := N \mathcal{E}(t) + \mathcal{F}_{1}(t) + \mathcal{F}_{2}(t) + \mathcal{F}_{3}(t) + N_{4} \mathcal{F}_{4}(t). \notag
		\end{equation}
		Using now the estimates for $\partial_{t} \mathcal{F}_{1}$, $\partial_{t} \mathcal{F}_{2}$, $\partial_{t} \mathcal{F}_{3}$ and $\partial_{t} \mathcal{F}_{4}$  
		from Equations (\ref{EQUATION_ABSCHAETZUNG_F1})--(\ref{EQUATION_ABSCHAETZUNG_F4}), we obtain
		\begin{equation}
			\begin{split}
				\partial_{t} \mathcal{L}(t) \leq& -C_{w_{t}} \int_{\Omega} w_{t}^{2} \mathrm{d} x
				-C_{\nabla w} \int_{\Omega} |\nabla w|^{2} \mathrm{d} x
				-C_{v_{t}} \int_{\Omega} |v_{t}|^{2} \mathrm{d} x
				-C_{\sqrt{S} \mathcal{D} v} \int_{\Omega} |\sqrt{S} \mathcal{D}' v|^{2} \mathrm{d} x \notag \\
				& -C_{\theta} \int_{\Omega} \theta^{2} \mathrm{d} x -C_{q} \int_{\Omega} |q|^{2} \mathrm{d} x
			\end{split}
			\notag
		\end{equation}
		with the constants
		\begin{equation}
			\begin{split}
				C_{w_{t}} &= d N - \tfrac{\rho_{1} + d}{2 \varepsilon_{1}} - \left(\tfrac{d}{2 \varepsilon_{2}} + \rho_{1}\right), \\
				C_{\nabla w} &= \left[K - \tfrac{\varepsilon_{2} C_{\mathcal{P}}}{2} (K + d)\right] - \tfrac{\rho_{3} K \varepsilon'_{3}}{2}, \\
				C_{v_{t}} &= \left[\rho_{2} - \tfrac{\varepsilon_{1}}{2} \left(\rho_{1} C_{\mathcal{P}}
				+ \tfrac{\gamma \tau_{0}}{\kappa}\right) \right]
				+ \left[\rho_{2} \gamma - \tfrac{\rho_{2} \kappa \varepsilon_{3}}{2}\right]
				- N_{4} \tfrac{\varepsilon'_{4} \tau_{0} \gamma}{2}, \\
				C_{\sqrt{S} \mathcal{D} v} &= \left[1 - \tfrac{\varepsilon_{1}}{2 C_{\mathcal{K}, 1}} \left(d C_{\mathcal{P}} + \tfrac{\gamma \delta}{\kappa} \right) \right] -
				\left(\tfrac{\rho_{3} 2 \|S\| \varepsilon_{3}'}{2} + \tfrac{\rho_{3} K C_{P} \varepsilon_{3}'}{2 C_{\mathcal{K}, 1}} \right), \\
				C_{\theta} &= \left(\rho_{3} \gamma - \tfrac{\rho_{3} (\|S\| + K) C_{\mathcal{B}}}{2 \varepsilon'_{3}} -
				\tfrac{\rho_{3} K C_{\mathcal{B}}}{2 \varepsilon'_{3}} \theta^{2}\right)
				+ N_{4} \left[\rho_{3} \kappa - \tfrac{\varepsilon_{4} \rho_{3} \delta C_{\mathcal{B}_{\mathrm{rot}}}}{2} \right], \\
				C_{q} &= \tau_{0} N - \tfrac{\gamma (\tau_{0} + \delta)}{2 \kappa \varepsilon_{1}} - \tfrac{\rho_{2} \kappa C'_{\mathcal{B}_{\mathrm{rot}}}}{2 \varepsilon_{3}} -
				N_{4} \left((1 + C'_{\mathcal{B}_{\mathrm{rot}}}) \tau_{0} \kappa + \tfrac{\rho_{3} \delta}{2 \varepsilon_{4}} + \tfrac{\tau_{0} \gamma}{2 \varepsilon'_{4}}\right).
			\end{split} \notag
		\end{equation}
		Now, we select $\varepsilon_{1}, \varepsilon_{2}, \varepsilon_{3}, \varepsilon_{4} > 0$ sufficiently small
		for all bracket terms in $C_{\sqrt{S} \mathcal{D} v}$, $C_{\nabla w}$, $C_{v_{t}}$ and $C_{\theta}$ to be positive.
		Next, we choose $\varepsilon_{3}' > 0$ so small
		that $C_{\nabla w}$ and $C_{\sqrt{S} \mathcal{D} v}$ become positive.
		Then we fix a sufficiently large $N_{4} > 0$ to assure for $C_{\theta} > 0$.
		We further pick a small $\varepsilon_{4}' > 0$ to make $C_{v_{t}}$ positive.
		Finally, we choose $N > 0$ to be sufficiently large
		to guarantee the positivity of $C_{w_{t}}$ and $C_{q}$.
		Therefore, we get
		\begin{equation}
			C_{\min} := \min\{C_{w_{t}}, C_{\nabla w}, C_{v_{t}}, C_{\sqrt{S} \mathcal{D} v}, C_{\theta}, C_{q}\} > 0. \notag
		\end{equation}
		Taking into account Young's inequality
		\begin{equation}
			|\nabla w + v|^{2} \leq \tfrac{1}{2}(|\nabla w|^{2} + |\nabla v|^{2}), \notag
		\end{equation}
		Korn's inequality immediately yields the estimate
		\begin{equation}
			\begin{split}
				|\nabla w| + |\sqrt{S} \mathcal{D} v|^{2}
				&\geq |\nabla w| + \tfrac{1}{2} C_{\mathcal{K}, 1} |v|^{2} + \tfrac{1}{2} |\sqrt{S} \mathcal{D} v|^{2}
				\geq \min\{2, C_{\mathcal{K}, 1}\} |\nabla w + v|^{2} + \tfrac{1}{2} |\sqrt{S} \mathcal{D} v|^{2} \\
				&\geq \min\{\tfrac{1}{2}, C_{\mathcal{K}, 1}\} (|\nabla w + v|^{2} + |\sqrt{S} \mathcal{D} v|^{2}).
			\end{split} \notag
		\end{equation}
		Hence, we get
		\begin{equation}
			\partial_{t} \mathcal{F}(t)
			\leq -2 \tfrac{\min\big\{1, \min\big\{\tfrac{1}{2}, C_{\mathcal{K}, 1}\big\}^{-1}\big\}}{
			\max\{1, \rho_{1}, \rho_{2}, \rho_{3}, \tau_{0}, K\}} \mathcal{E}(t) =: C \mathcal{E}(t).
			\notag
		\end{equation}
		On the other hand, we can estimate
		\begin{equation}
			\begin{split}
				|\mathcal{F}_{1} + \mathcal{F}_{2} + \mathcal{F}_{3} + &N_{4} \mathcal{F}_{4}|(t)
				\leq
				\tfrac{1}{2} \int_{\Omega} \bigg(\rho_{1} (w_{t}^{2} + |u|^{2}) + \rho_{2} (|v_{t}|^{2} + |v|^{2})
				+ \tfrac{\gamma \tau_{0}}{\kappa} (|v|^{2} + |q|^{2}) + \\
				&\rho_{1} (w_{t}^{2} + w^{2}) +
				\rho_{2} \rho_{3} (|\mathcal{B}_{\mathrm{rot}} \theta|^{2} + |v_{t}|^{2})
				+ \tau_{0} \rho_{3} (|q|^{2} + |\mathcal{B}_{\mathrm{rot}} \theta|^{2})\bigg) \mathrm{d} x \\
				\leq& \tfrac{1}{2} \bigg(2\rho_{1} \|w_{t}\|_{L^{2}(\Omega)}^{2} +
				\rho_{1} \|w\|_{H^{1}(\Omega)}^{2} + \rho_{2} \|v_{t}\|_{(L^{2}(\Omega)^{2})}^{2}
				+ (\rho_{2} + \tfrac{\gamma \tau_{0}}{\kappa}) \|v\|_{(H^{1}(\Omega))^{2}} \\
				&+ C_{\mathcal{B}_{\mathrm{rot}}} (\rho_{2} \rho_{3} + \tau_{0} \rho_{3}) \|\theta\|_{L^{2}(\Omega)}^{2}
				+ (\tfrac{\gamma \tau_{0}}{\kappa} + \tau_{0} \rho_{3}) \|q\|_{(L^{2}(\Omega))^{2}}^{2}\bigg) \\
				\leq& \tfrac{1}{2} \bigg(2\rho_{1} \|w_{t}\|_{L^{2}(\Omega)}^{2} + \rho_{2} \|v_{t}\|_{(L^{2}(\Omega)^{2})}^{2} + \\
				&\tfrac{\max\{\rho_{1}, (\rho_{2} + \tfrac{\gamma \tau_{0}}{\kappa})\}}{C_{\mathcal{K}}} (
				K \|\nabla w + v\|_{(L^{2}(\Omega))^{2}}^{2} + \|\sqrt{S} \mathcal{D} v\|_{(L^{2}(\Omega))^{2}}^{2})
				\\
				&+ C_{\mathcal{B}_{\mathrm{rot}}} (\rho_{2} \rho_{3} + \tau_{0} \rho_{3}) \|\theta\|_{L^{2}(\Omega)}^{2}
				+ (\tfrac{\gamma \tau_{0}}{\kappa} + \tau_{0} \rho_{3}) \|q\|_{(L^{2}(\Omega))^{2}}^{2}\bigg)
				\leq \hat{C} \mathcal{E}(t).
			\end{split} \notag
		\end{equation}
		Letting now
		$\alpha_{1} := N - \tfrac{\max\{\rho_{1}, \rho_{2}, C_{\mathcal{K}}^{-1}\}}{\min\{\rho_{1}, \rho_{2}, \rho_{3}\}}$ and
		$\alpha_{2} := N + \tfrac{\max\{\rho_{1}, \rho_{2}, C_{\mathcal{K}}^{-1}\}}{\min\{\rho_{1}, \rho_{2}, \rho_{3}\}}$,
		we obtain the following equivalence between $\mathcal{E}$ and $\mathcal{F}$
		\begin{equation}
			\alpha_{1} \mathcal{E}(t) \leq \mathcal{F}(t) \leq \alpha_{2} \mathcal{E}(t) \text{ for } t \geq 0. \notag
		\end{equation}
		If necessary, we increase the constant $N$ to assure for the positivity of $\alpha_{1}$.
		Thus, both $C$, $\alpha_{1}$ and $\alpha_{2}$ are positive.
		Exploiting Gronwall's inequality, we obtain the following estimate for $\mathcal{E}$
		\begin{equation}
			\mathcal{E}(t) \leq \tfrac{1}{\alpha_{1}} \mathcal{F}(t)
			\leq \tfrac{1}{\alpha_{1}} \mathcal{E}(0) e^{-\tfrac{C}{\alpha_{2}} t}
			=: C \mathcal{E}(0) e^{-2\alpha t} \text{ for } t \geq 0 \notag
		\end{equation}
		meaning an exponential decay of $\mathcal{E}$.
	\end{proof}

\begin{appendices}
\section{The divergence problem and the Bogowski\u i operator}
	In various applications of partial differential equations, e.g., when studying Navier-Stokes equations,
	there arises a so-called ``divergence problem'': 
	For a given function $f$, determine a vector field $u$ such that its divergence coincides with $f$.
	We refer to \cite{GeiHeHie2006} for a rather general solution of this problem in bounded domains.
	It has namely been shown that the solution map $\mathcal{B} \colon f \mapsto u$, called the Bogowski\u i-operator, is a bounded linear operator between
	$W^{s, p}_{0}(\Omega)$ and $W^{s + 1, p}_{0}(\Omega)$ for $p \in (0, \infty)$, $s \in \big(-2 + \tfrac{1}{p}, \infty\big)$.
        
	For our application, we want to additionally guarantee that the solution $u$ is irrotational.
	To this end, we exploit the following result from \cite{Ir2006}.
	\begin{theorem}
		\label{SATZ_GRADIENT_DIVERGENZ_ROTATIONSFREIHEIT}
		Let $\Omega \subset \mathbb{R}^{n}$ be a domain with a smooth boundary
		and let $\nu \colon \Omega \to \mathbb{R}^{n}$ denote the outer unit normal vector on $\partial \Omega$.
		There exists then a function $u \in H^{1}(\Omega, \mathbb{R}^{n})$ satisfying $\nu \otimes u = u \otimes \nu$ on $\partial \Omega$ and
		\begin{equation}
			\|\nabla u\|_{L^{2}(\Omega)}^{2} = \|\mathrm{div}\, u\|_{L^{2}(\Omega)}^{2} + \tfrac{1}{2} \|\nabla u - (\nabla u)'\|_{L^{2}(\Omega)}^{2}
			+ (n - 1) \int_{\partial \Omega} |u|^{2} H_{n} \mathrm{d}S,
			\label{GLEICHUNG_NORM_DES_GRADIENTEN_GEGEN_DIVERGENZ_UND_ROTATION}
		\end{equation}
		where $H_{n} \colon \partial \Omega \to \mathbb{R}$, $x \mapsto H_{n}(x)$ denotes the mean curvature of $\partial \Omega$ with respect to the outer normal vector.
		In $n = 2, 3$, Equation (\ref{GLEICHUNG_NORM_DES_GRADIENTEN_GEGEN_DIVERGENZ_UND_ROTATION}) reduces to
		\begin{equation}
			\|\nabla u\|_{L^{2}(\Omega)}^{2} = \|\mathrm{div}\, u\|_{L^{2}(\Omega)}^{2} + \|\mathrm{rot} u\|_{L^{2}(\Omega)}^{2}
			+ (n - 1) \int_{\partial \Omega} |u|^{2} H_{n} \mathrm{d}S,
			\label{GLEICHUNG_NORM_DES_GRADIENTEN_GEGEN_DIVERGENZ_UND_ROTATION_2D_UND_3D}
		\end{equation}
		where
		\begin{equation}
			\mathrm{rot} u =
			\begin{pmatrix}
				\partial_{x_{2}} u_{3} - \partial_{x_{3}} u_{2} \\
				\partial_{x_{3}} u_{1} - \partial_{x_{1}} u_{3} \\
				\partial_{x_{1}} u_{2} - \partial_{x_{2}} u_{1}
			\end{pmatrix} \text{ for } n = 3 \text{ and }
			\mathrm{rot} u = \partial_{x_{1}} u_{2} - \partial_{x_{2}} u_{1} \text{ for } n = 2. \notag
		\end{equation}
		For $u \in H^{1}_{0}(\Omega, \mathbb{R}^{n})$, the second term in
		(\ref{GLEICHUNG_NORM_DES_GRADIENTEN_GEGEN_DIVERGENZ_UND_ROTATION}) and (\ref{GLEICHUNG_NORM_DES_GRADIENTEN_GEGEN_DIVERGENZ_UND_ROTATION_2D_UND_3D})
		vanishes and no assumptions on $\partial \Omega$ are required:
		\begin{equation}
			\|\nabla u\|_{L^{2}(\Omega)}^{2} = \|\mathrm{div}\, u\|_{L^{2}(\Omega)}^{2} + \|\nabla u - (\nabla u)'\|_{L^{2}(\Omega)}^{2}.
			\label{GLEICHUNG_NORM_DES_GRADIENTEN_GEGEN_DIVERGENZ_UND_ROTATION_2D_UND_3D_OHNE_RANDTERME}
		\end{equation}
	\end{theorem}
	
	In the following, we assume $n = 2$.
	We define the space
	\begin{equation}
		H^{1}_{0, \mathrm{rot}}(\Omega) = \big\{u \in (H^{1}_{0}(\Omega))^{2} \,|\, \nabla u = (\nabla u)'\big\} =
		\big\{u \in (H^{1}_{0}(\Omega))^{2} \,|\, \mathrm{rot} u = 0 \big\}
		\notag
	\end{equation}
	equipped with the standard inner product of $(H^{1}_{0}(\Omega))^{2}$.
	Since $H^{1}_{0, \mathrm{rot}}(\Omega)$ is a closed subspace of $(H^{1}_{0}(\Omega))^{2}$, $H^{1}_{0, \mathrm{rot}}(\Omega)$ is a Hilbert space.
	We prove the following theorem.
	\begin{theorem}
		\label{SATZ_BOGOWSKII_STETIGKEIT}
		The mapping
		\begin{equation}
			\mathrm{div}\, \colon H^{1}_{0, \mathrm{rot}}(\Omega) \to L^{2}(\Omega)/\{1\} \notag
		\end{equation}
		is an isomorphism with an inverse $\mathrm{div}\,^{-1} = \mathcal{B}_{\mathrm{rot}\,}$
		\begin{equation}
			\mathcal{B}_{\mathrm{rot}} \colon L^{2}(\Omega)/\{1\} \to H^{1}_{0, \mathrm{rot}}(\Omega) \notag
		\end{equation}
		in the sense
		\begin{equation}
			\mathrm{div}\, \mathcal{B}_{\mathrm{rot}} = \mathrm{id}_{L^{2}(\Omega)/\{1\}} \text{ and }
			\mathcal{B}_{\mathrm{rot}} \mathrm{div}\, = \mathrm{id}_{H^{1}_{0, \mathrm{rot}}(\Omega)}. \notag
		\end{equation}
		Furthermore, the exists $C_{\mathcal{B}} > 0$ such that
		\begin{equation}
			\|\mathcal{B}_{\mathrm{rot}} f\|_{(H^{1}(\Omega))^{2}} \leq C_{\mathcal{B}_{\mathrm{rot}}} \|f\|_{L^{2}(\Omega)}
			\notag
		\end{equation}
		holds true for all $f \in L^{2}_{\ast}(\Omega)$.
        \end{theorem}
	\begin{proof}
		The linearity of $\mathrm{div}$ is obvious
		For each $u \in H^{1}_{0, \mathrm{rot}}(\Omega)$, we have $\mathrm{div}\, u \in L^{2}(\Omega)$ and thus
		\begin{equation}
			\int_{\Omega} \mathrm{div}\, u \mathrm{d}x = \int_{\Gamma} u \cdot \nu \mathrm{d}\Gamma = 0, \notag
		\end{equation}
		meaning $\mathrm{div}\, u \in L^{2}(\Omega)/\{1\}$.
		The continuity is also trivial since
		\begin{equation}
			\|\mathrm{div}\, u\|_{L^{2}(\Omega)} \leq \sqrt{2} \|\nabla u\|_{L^{2}(\Omega)} \leq \sqrt{2} \|u\|_{H^{1}(\Omega)}. \notag
		\end{equation}
		The operator $\mathrm{div}\,$ is injective. Indeed, let $u_{1}, u_{2} \in H^{1}_{0, \mathrm{rot}}(\Omega)$.
		Let $\mathrm{div}\, u_{1} = \mathrm{div}\, u_{2}$. Then, using Poincar\'{e} inequality,
		\begin{equation}
			0 = \|\mathrm{div}\, u_{1} - \mathrm{div}\, u_{2}\|_{L^{2}(\Omega)} \geq \|\nabla u_{1} - \nabla u_{2}\|_{L^{2}(\Omega)}
			\geq \tfrac{1}{C_{\mathcal{P}}} \|u_{1} - u_{2}\|_{L^{2}(\Omega)}, \notag
		\end{equation}
		i.e., $u_{1} = u_{2}$.
	
		To explicitely construct the operator $\mathcal{B}_{\mathrm{rot}}$, we follow the variational approach.
		For $f, g \in L^{2}(\Omega)/\{1\}$, we consider a boundary value problem for $\varphi, \psi \in H^{1}(\Omega)/\{1\}$:
		\begin{equation}
			\begin{split}
				-\mathrm{div}\, (\nabla \varphi + \mathrm{rot}' \psi) &= f \text{ in } \Omega, \\
				-\mathrm{rot} (\nabla \varphi + \mathrm{rot}' \psi) &= g \text{ in } \Omega, \\
				\nu \cdot (\nabla \varphi + \mathrm{rot}' \psi) &= 0 \text{ on } \Gamma, \\
				\nu^{\perp} \cdot (\nabla \varphi + \mathrm{rot}' \psi) &= 0 \text{ on } \Gamma, \\
			\end{split} \label{GLEICHUNG_HELMHOLZ}
		\end{equation}
		where $\nu^{\perp} := (\nu_{2}, -\nu_{1})'$, $\mathrm{rot}' := (\partial_{x_{2}}, -\partial_{x_{1}})'$.
		We multiply the equations with $\tilde{\varphi}, \tilde{\psi} \in H^{1}(\Omega)/\{1\}$,
		sum up the resulting identities, take into account the boundary conditions and apply a partial integration to find
		\begin{equation}
			-\int_{\Omega} \mathrm{div}\, (\nabla \varphi + \mathrm{rot}' \psi) \tilde{\varphi} \mathrm{d}x
			-\int_{\Omega} \mathrm{rot} (\nabla \varphi + \mathrm{rot}' \psi) \tilde{\psi} \mathrm{d}x =
			\int_{\Omega} (\nabla \varphi + \mathrm{rot}' \psi) \cdot (\nabla \tilde{\varphi} + \mathrm{rot} \tilde{\psi}) \mathrm{d}x \notag
		\end{equation}
		This lead to the following operator equation
		\begin{equation}
			\mathcal{A} (\varphi, \psi)' = (f, g)',
			\label{GLEICHUNG_BOGOWSKI_OPERATORFORMULIERUNG_ORIGINAL_ANHANG}
		\end{equation}
		where
		\begin{equation}
			\mathcal{A} \colon D(\mathcal{A}) \subset \mathcal{H} \to \mathcal{H}, \quad
			(\varphi, \psi)' \mapsto
			\begin{pmatrix}
				-\mathrm{div}\, (\nabla \varphi + \mathrm{rot}' \psi) \\
				-\mathrm{rot} (\nabla \varphi + \mathrm{rot}' \psi)
			\end{pmatrix} \notag
		\end{equation}
		and
		\begin{equation}
			D(\mathcal{A}) = \Big\{(\varphi, \psi)' \in \mathcal{V} \,\big|\,
			\exists (f_{1}, f_{2})' \in \mathcal{H} \; \forall (\tilde{\varphi}, \tilde{\psi})' \in \mathcal{V}:
			B(\varphi, \psi; \tilde{\varphi}, \tilde{\psi}) = \int_{\Omega} f_{1} \tilde{\varphi} + f_{2} \tilde{\psi} \mathrm{d} x \Big\}
			\notag
		\end{equation}
		with the bilinear form
		\begin{equation}
			B \colon \mathcal{V} \times \mathcal{V} \to \mathbb{R}, \quad
			(\phi, \psi, \tilde{\phi}, \tilde{\psi})' \mapsto
			\int_{\Omega} (\nabla \varphi + \mathrm{rot}' \psi) \cdot (\nabla \tilde{\varphi} + \mathrm{rot} \tilde{\psi}) \mathrm{d}x.
			\notag
		\end{equation}
		Here, we introduced the Hilbert spaces
		\begin{align}
			\mathcal{H} := (L^{2}(\Omega)/\{1\}) \times (L^{2}(\Omega)/\{1\}), \quad
			\mathcal{V} := (H^{1}(\Omega)/\{1\}) \times (H^{1}(\Omega)/\{1\}) \notag
		\end{align}
		equipped with the standard inner products of $L^{2}(\Omega) \times L^{2}(\Omega)$ and $H^{1}(\Omega) \times H^{1}(\Omega)$, respectively.
		Since $\mathcal{A}$ has a nontrivial kernel, we consider the operator given as its restriction onto the closed subspace
		\begin{equation}
			\tilde{\mathcal{V}} = \{(\varphi, \psi)' \in \mathcal{V} \,|\, \forall (\tilde{\varphi}, \tilde{\psi})' \in \mathcal{V}:
			\int_{\Omega} \nabla \varphi \cdot \mathrm{rot}' \tilde{\psi} \mathrm{d} x =
			\int_{\Omega} \nabla \tilde{\varphi} \cdot \mathrm{rot}' \psi \mathrm{d} x = 0 \}\notag
		\end{equation}
		of $\mathcal{V}$ and denote it as
		\begin{equation}
			\tilde{\mathcal{A}} \colon D(\tilde{\mathcal{A}}) := D(\mathcal{A}) \cap \tilde{\mathcal{V}} \subset \mathcal{H} \to \mathcal{H}. \notag
		\end{equation}
		Equation (\ref{GLEICHUNG_BOGOWSKI_OPERATORFORMULIERUNG_ORIGINAL_ANHANG}) reduces then to
		\begin{equation}
			\tilde{\mathcal{A}} (\varphi, \psi)' = (f, g)'.
			\label{GLEICHUNG_BOGOWSKI_OPERATORFORMULIERUNG_ANHANG}
		\end{equation}
		We multiply Equation (\ref{GLEICHUNG_BOGOWSKI_OPERATORFORMULIERUNG_ANHANG}) scalar in $\mathcal{H}$ with $(\tilde{\varphi}, \tilde{\psi})' \in \tilde{\mathcal{V}}$
		to find after a partial integration the weak formulation of (\ref{GLEICHUNG_BOGOWSKI_OPERATORFORMULIERUNG_ANHANG}):
		Determine an element $(\varphi, \psi)' \in \tilde{\mathcal{V}}$ such that
		\begin{equation}
			B(\varphi, \psi; \hat{\varphi}, \hat{\psi}) = F(\hat{\varphi}, \hat{\psi})
			\text{ for all } (\hat{\varphi}, \hat{\psi})' \in \tilde{\mathcal{V}},
			\label{GLEICHUNG_BOGOWSKI_SCHWACHE_FORMULIERUNG_ANHANG}
		\end{equation}
		where
		\begin{equation}
			\begin{split}
				B &\colon \tilde{\mathcal{V}} \times \tilde{\mathcal{V}} \to \mathbb{R}, \quad
				(\phi, \psi, \hat{\phi}, \hat{\psi})' \mapsto
				\int_{\Omega} (\nabla \varphi + \mathrm{rot}' \psi) \cdot (\nabla \hat{\varphi} + \mathrm{rot} \hat{\psi}) \mathrm{d}x, \\
				F &\colon \tilde{\mathcal{V}} \to \mathbb{R}, \quad
				(\hat{\phi}, \hat{\psi})' \mapsto \int_{\Omega} \hat{\varphi} f \mathrm{d} x + \int_{\Omega} \hat{\psi} g \mathrm{d} x.
			\end{split} \notag
		\end{equation}
		The bilinear form $B$ and the linear functional $F$ are continuous on $\tilde{\mathcal{V}} \times \tilde{\mathcal{V}}$ and $\tilde{\mathcal{V}}$, respectively.
		The bilinear form $B$ is symmetrical.
		By the virtue of second \textsc{Poincar\'{e}}'s inequality, we obtain
		\begin{align}
			B(\varphi, \psi) &=
			\|\nabla \varphi\|_{L^{2}(\Omega)}^{2} + 2 \langle \nabla \varphi, \mathrm{rot}' \psi\rangle + \|\mathrm{rot}' \psi\|_{L^{2}(\Omega)}^{2} \notag \\
			&= \|\nabla \varphi\|_{L^{2}(\Omega)}^{2} + \|\mathrm{rot}' \psi\|_{L^{2}(\Omega)}^{2} =
			\|\nabla \varphi\|_{L^{2}(\Omega)}^{2} + \|\nabla \psi\|_{L^{2}(\Omega)}^{2} \notag \\
			&\geq \tfrac{1}{2}(1 + \tfrac{1}{C_{\mathcal{P}}}) (\|\varphi\|_{H^{1}(\Omega)}^{2} + \|\psi\|_{H^{1}(\Omega)}^{2})
			= \tfrac{1}{2}(1 + \tfrac{1}{C_{\mathcal{P}}}) \|(\varphi, \psi)'\|_{\mathcal{V}}^{2} =: b \|(\varphi, \psi)'\|_{\tilde{\mathcal{V}}}^{2} , \notag
		\end{align}
		i.e., $B$ is coercive.
		The lemma of Lax \& Milgram yields the existence of a unique solution $(\varphi, \psi)' \in \tilde{\mathcal{V}}$ to Equation (\ref{GLEICHUNG_BOGOWSKI_SCHWACHE_FORMULIERUNG_ANHANG}).
		There further holds
		\begin{equation}
			b \|(\varphi, \psi)'\|_{\tilde{\mathcal{V}}}^{2} \leq B(\varphi, \psi)
			\leq \tfrac{b}{2} \|(\varphi, \psi)'\|_{\mathcal{H}}^{2} + \tfrac{1}{2b} \|(f, g)'\|_{\mathcal{H}}^{2}
			\leq \tfrac{b}{2} \|(\varphi, \psi)'\|_{\tilde{\mathcal{V}}}^{2} + \tfrac{1}{2b} \|(f, g)'\|_{\mathcal{H}}^{2}, \notag
		\end{equation}
		i.e.,
		\begin{equation}
			\|(\varphi, \psi)'\|_{\tilde{\mathcal{V}}}^{2} \leq \tfrac{1}{b} \|(f, g)'\|_{\mathcal{H}}^{2}. \notag
		\end{equation}
		Exploiting the trivial identities
		\begin{equation}
			\mathrm{div}\, \, \mathrm{rot}' \varphi = 0, \quad
			\mathrm{rot} \, \nabla \varphi = 0 \text{, etc., in } (\mathcal{C}_{0}^{\infty}(\Omega))' \notag
		\end{equation}
		and the definition of $\mathcal{V}$, we find
		\begin{equation}
			\int_{\Gamma} \nu \varphi \cdot \hat{\mathrm{rot}}' \psi \mathrm{d} x
			= \int_{\Omega} \nabla \varphi \cdot \mathrm{rot}' \hat{\psi} \mathrm{d} x = 0, \quad
			\int_{\Gamma} \nu^{\perp} \psi \cdot \nabla \hat{\varphi} \mathrm{d} \Gamma =
			\int_{\Gamma} \mathrm{rot}' \psi \cdot \nabla \hat{\varphi} \mathrm{d}x = 0 \text{, etc.} \notag
		\end{equation}
		for all $(\varphi, \psi)' \in \tilde{\mathcal{V}}$ and $(\hat{\varphi}, \hat{\psi})' \in \mathcal{V}$.
		Hence,
		\begin{equation}
			\begin{split}
				-\int_{\Omega} \mathrm{div}\, (\nabla \varphi + &\mathrm{rot}' \psi) \hat{\varphi}
				+ \mathrm{rot} (\nabla \varphi + \mathrm{rot}' \psi) \hat{\psi} \mathrm{d} x \\
				= &B(\varphi, \psi; \hat{\varphi}, \hat{\psi})
				-\int_{\Gamma} \nu \cdot (\nabla \varphi + \mathrm{rot}' \psi) \hat{\varphi} +
				\nu^{\perp} \cdot (\nabla \varphi + \mathrm{rot}' \psi) \hat{\varphi} \mathrm{d}\Gamma
			\end{split} \notag
		\end{equation}
		holds true for all $(\hat{\varphi}, \hat{\psi})' \in \mathcal{V}$
		and, in particular, the solution $(\varphi, \psi)' \in \tilde{\mathcal{V}}$ of (\ref{GLEICHUNG_BOGOWSKI_SCHWACHE_FORMULIERUNG_ANHANG}).
		Therefore, $(\varphi, \psi)' \in D(\tilde{\mathcal{A}})$.
		Thus, we have shown that $\tilde{\mathcal{A}}$ is invertible
		and its inverse $\tilde{\mathcal{A}}^{-1} \colon \mathcal{H} \to D(\tilde{\mathcal{A}})$ is continuous:
		\begin{equation}
			\|\tilde{\mathcal{A}}^{-1} (f, g)'\|_{\mathcal{V}} \leq \tfrac{1}{b} \|(f, g)'\|_{\mathcal{H}}^{2} \notag
		\end{equation}
		
		Let $f \in L^{2}(\Omega)/\{1\}$. We define $(\phi, \psi) := \tilde{\mathcal{A}}^{-1} (f, 0)'$, $u := \nabla \varphi + \mathrm{rot}' \psi$
		and obtain by construction
		\begin{equation}
			\begin{split}
				\mathrm{div}\, u &= \triangle \varphi = f \text{ in } \Omega, \\
				\mathrm{rot} u &= \mathrm{rot} 0 = 0 \text{ in } \Omega, \\
				u &= \nabla \varphi + \mathrm{rot} \psi = 0 \text{ on } \Gamma,
			\end{split}
		\end{equation}
		i.e., $u \in H^{1}_{\mathrm{rot}}(\Omega)$ with $\mathrm{div}\, u = f$.
		Thus, there exists a continuous inverse
		\begin{equation}
			\mathcal{B}_{\mathrm{rot}} \colon L^{2}(\Omega)/\{1\} \to H^{1}_{0, \mathrm{rot}}(\Omega), \quad
			f \mapsto u \notag
		\end{equation}
		of $\mathrm{div}\,$ such that
		\begin{equation}
			\begin{split}
				\|\mathcal{B}_{\mathrm{rot}} f\|_{(H^{1}(\Omega))^{2}}
				&= \|\mathcal{B}_{\mathrm{rot}} f\|_{(L^{2}(\Omega))^{2}}^{2} + \|\nabla \mathcal{B}_{\mathrm{rot}} f\|_{(L^{2}(\Omega))^{2 \times 2}}^{2} \\
				&= \|\nabla \varphi + \mathrm{rot}' \psi \|_{(L^{2}(\Omega))^{2}}^{2} + \|\mathrm{div}\, \mathcal{B}_{\mathrm{rot}} f\|_{L^{2}(\Omega)}^{2} \\
				&\leq 2\|\nabla \varphi\|_{(L^{2}(\Omega))^{2}}^{2} + 2\|\mathrm{rot}' \psi\|_{(L^{2}(\Omega))^{2}}^{2}
				+ \|f\|_{L^{2}(\Omega)}^{2} \\
				&\leq (\tfrac{2}{b} + 1) \|f\|_{L^{2}(\Omega)}^{2} =: C_{\mathcal{B}_{\mathrm{rot}}} \|f\|_{L^{2}(\Omega)}.
			\end{split}
			\notag
		\end{equation}
		This finishes the proof.
	\end{proof}
	
	\begin{corollary}
		The operator $\mathcal{B}_{\mathrm{rot}}$ can be extended to a linear continuous operator 
		\begin{equation}
			\mathcal{B}_{\mathrm{rot}} \colon (H^{1}(\Omega))' \to (L^{2}(\Omega))^{2}. \notag
		\end{equation}
		(Cp. also \cite{BoSo1990, GeiHeHie2006} for the rotational case.)
	\end{corollary}

	\begin{proof}
		Due to the coercivity of the bilinear form $B$, the operator $\tilde{\mathcal{A}}$ defined in the proof of Theorem \ref{SATZ_BOGOWSKII_STETIGKEIT} strictly positive.
		According to \cite[Section 3.4]{TuWei2009}, it is possible to define square roots
		\begin{equation}
			\tilde{\mathcal{A}}^{-1/2} \in L(\mathcal{H}, \mathcal{H}) \text{ and } \tilde{\mathcal{A}}^{1/2} \colon D(\tilde{\mathcal{A}}^{1/2}) := \mathrm{im} \, \tilde{\mathcal{A}}^{-1/2} \to \mathcal{H} \notag
		\end{equation}
		of $\tilde{\mathcal{A}}^{-1}$ and $\tilde{\mathcal{A}}$, respectively.
		Further, there exists a continuous continuation of $\tilde{\mathcal{A}}^{-1}$
		\begin{equation}
			\tilde{\mathcal{A}}^{-1} \in L(D(\tilde{\mathcal{A}}^{-1/2}), D(\tilde{\mathcal{A}}^{1/2})), \notag
		\end{equation}
		where $D(\tilde{\mathcal{A}}^{-1/2}) = D(\tilde{\mathcal{A}}^{1/2})'$.
		Hence,
		\begin{equation}
			\tilde{\mathcal{B}}_{\mathrm{rot}} \colon D(\tilde{\mathcal{A}}^{-1/2}) \to (L^{2}(\Omega))^{2}, \quad
			f \mapsto \nabla \varphi + \mathrm{rot}' \psi
			\text{ with } (\varphi, \psi)' := \tilde{\mathcal{A}}^{-1}(f, 0)' \in \tilde{\mathcal{V}} \notag
		\end{equation}
		represents a continuous continuation of $\mathcal{B}_{\mathrm{rot}}$ onto $D(\tilde{\mathcal{A}}^{-1/2})$.
		Since $(H^{1}(\Omega))' \subset D(\tilde{\mathcal{A}}^{-1/2})$
		and the norms of $(H^{1}(\Omega))'$ und $D(\tilde{\mathcal{A}}^{-1/2})$ are equivalent, the claim follows.
	\end{proof}
	
	Let us now consider a vector field $u \in (H^{1}(\Omega))^{2}$ with $u \cdot \nu = 0$ on $\Gamma$.
	Unfortunately, the identity
	\begin{equation}
		\mathcal{B}_{\mathrm{rot}} \mathrm{div}\, u = u \notag
	\end{equation}
	does not hold in general since $u$ is not necessarily an element of $H^{1}_{0, \mathrm{rot}}(\Omega)$.
	Nevertheless, the following estimate holds true.
	\begin{theorem}
		\label{SATZ_BOGOWSKII_ABSCHAETZUNG}
		Let $u \in H^{1}(\Omega)$ satisfy $u \cdot \nu = 0$ on $\Gamma$.
		There exists then a constant $C'_{\mathcal{B}_{\mathrm{rot}}} > 0$ such that
		\begin{equation}
			\|\mathcal{B}_{\mathrm{rot}} \mathrm{div}\, u\|_{L^{2}(\Omega)} \leq C'_{\mathcal{B}_{\mathrm{rot}}} \|u\|_{(L^{2}(\Omega))^{2}} \notag
		\end{equation}
		for any $u \in (H^{1}(\Omega))^{2}$.
	\end{theorem}
	
	\begin{proof}
		We can estimate
		\begin{align}
			\|\mathcal{B}_{\mathrm{rot}} \mathrm{div}\, u\|_{(L^{2}(\Omega))^{2}} \leq C_{\mathcal{B}_{\mathrm{rot}}}
			\|\mathrm{div}\, u\|_{H^{-1}(\Omega)}. \notag
		\end{align}
		Further, we find
		\begin{equation}
			\int_{\Omega} \mathrm{div}\, u f \mathrm{d} x = -\int_{\Omega} u \nabla f \mathrm{d} x + \int_{\partial \Omega} u \cdot \nu f \mathrm{d} \Gamma
			= -\int_{\Omega} u \nabla f \mathrm{d} x
		\end{equation}
		for all $f \in H^{1}(\Omega)$ and therefore
		\begin{align}
			\|\mathrm{div}\, u\|_{H^{-1}(\Omega)}
			&= \sup_{\|f\|_{H^{1}(\Omega)} = 1}
			\left|\int_{\Omega} \mathrm{div}\, u f \mathrm{d} x\right|
			= \sup_{\|f\|_{H^{1}(\Omega)} = 1}
			\left|\int_{\Omega} u \nabla f \mathrm{d} x\right| \notag \\
			&\leq \sup_{\|f\|_{H^{1}(\Omega)} = 1} \|u\|_{(L^{2}(\Omega))^{2}} \|f\|_{H^{1}(\Omega)}
			= \|u\|_{(L^{2}(\Omega))^{2}}. \notag
		\end{align}
		This yields
		\begin{equation}
			\|\mathcal{B}_{\mathrm{rot}} \mathrm{div}\, u\|_{L^{p}(\Omega)} \leq C'_{\mathcal{B}_{\mathrm{rot}}} \|u\|_{L^{p}(\Omega)} \text{ for all } u \in H^{1}(\Omega)
			\label{GLEICHUNG_BOGOWSKII_DIVERGENZABSCHAETZUNG}
		\end{equation}
		with $C'_{\mathcal{B}_{\mathrm{rot}}} = C_{\mathcal{B}_{\mathrm{rot}}}$.
	\end{proof}
\end{appendices}

\section*{Acknowledgment}
The present work is dedicated to the honorable Mr. Urs Schaubhut, J.D. (Konstanz, Germany) in deep gratitude for his invaluable support in the author's struggle for justice.

\end{document}